\documentclass[10pt]{article}
\usepackage[latin1]{inputenc}
\usepackage{epsfig}
\usepackage{color}
\usepackage[british,english]{babel}
\usepackage{amsthm}
\usepackage{amsmath}
\usepackage{amsfonts}
\usepackage{amssymb}
\usepackage{graphicx}
\usepackage{multirow}

\usepackage{fancyhdr}

\newtheorem{corollary}{Corollary}[section]
\newtheorem{definition}[corollary]{Definition}

\newtheorem{lemma}[corollary]{Lemma}
\newtheorem{proposition}[corollary]{Proposition}
\newtheorem{remark}[corollary]{Remark}
\newtheorem{theorem}[corollary]{Theorem}

\date{}
\begin{document}
\title{Sharp pressure estimates for the Navier-Stokes system\\
 in thin porous media}\maketitle

\vskip-30pt
 \centerline{Mar\'ia ANGUIANO\footnote{Departamento de An\'alisis Matem\'atico. Facultad de Matem\'aticas.
Universidad de Sevilla, 41012 Sevilla (Spain)
anguiano@us.es} and  Francisco J. SU\'AREZ-GRAU\footnote{Departamento de Ecuaciones Diferenciales y An\'alisis Num\'erico. Facultad de Matem\'aticas. Universidad de Sevilla, 41012 Sevilla (Spain) fjsgrau@us.es}}

 \renewcommand{\abstractname} {\bf Abstract}
\begin{abstract} 
 A relevant problem for applications is to model the behavior of Newtonian fluids through thin porous media, which is a domain with small thickness $\epsilon$ and perforated by periodically distributed cylinders of size and period $\epsilon^\delta$, with $\delta>0$.   Depending on  the relation between thickness and the size of the cylinders, it was introduced in (Fabricius {\it et al.}, Transp. Porous Media, 115, 473-493, 2016), (Anguiano and Su\'arez-Grau, Z. Angew. Math. Phys., 68:45, 2017) and (Anguiano and Su\'arez-Grau, Mediterr. J. Math., 15:45, 2018) that there exist three regimes depending on the value of $\delta$: $\delta\in (0,1)$, $\delta=1$ and $\delta>1$.  In each regime, the asymptotic behavior of the fluid is governed by a lower-dimensional Darcy's law.\\
 
 In previous studies, the Reynolds number is considered to be of order one and so, the question that arises is for what range of values of the Reynolds number the lower-dimensional Darcy laws are still valid in each regime, which represents the main the goal of this paper. In this sense, considering a fluid governed by the Navier-Stokes system and assuming the Reynolds number written in terms of the thickness $\epsilon$, we prove that, for each regime, there exists a critical Reynolds number $Re_c$ such that for every Reynolds number $Re$ with order smaller or equal than $Re_c$, the lower-dimensional Darcy law is still valid. On the contrary, for Reynolds numbers $Re$ greater than $Re_c$, the inertial term of the Navier-Stokes system has to be taken into account in the asymptotic behavior and so, the Darcy law is not valid.

\end{abstract}

 {\small \bf AMS classification numbers: } 35Q30, 35B27, 74K35, 76S05, 76F06.
 
\section {Introduction}\label{S1} 
\begin{figure}[h]
\centering\includegraphics[width=5cm]{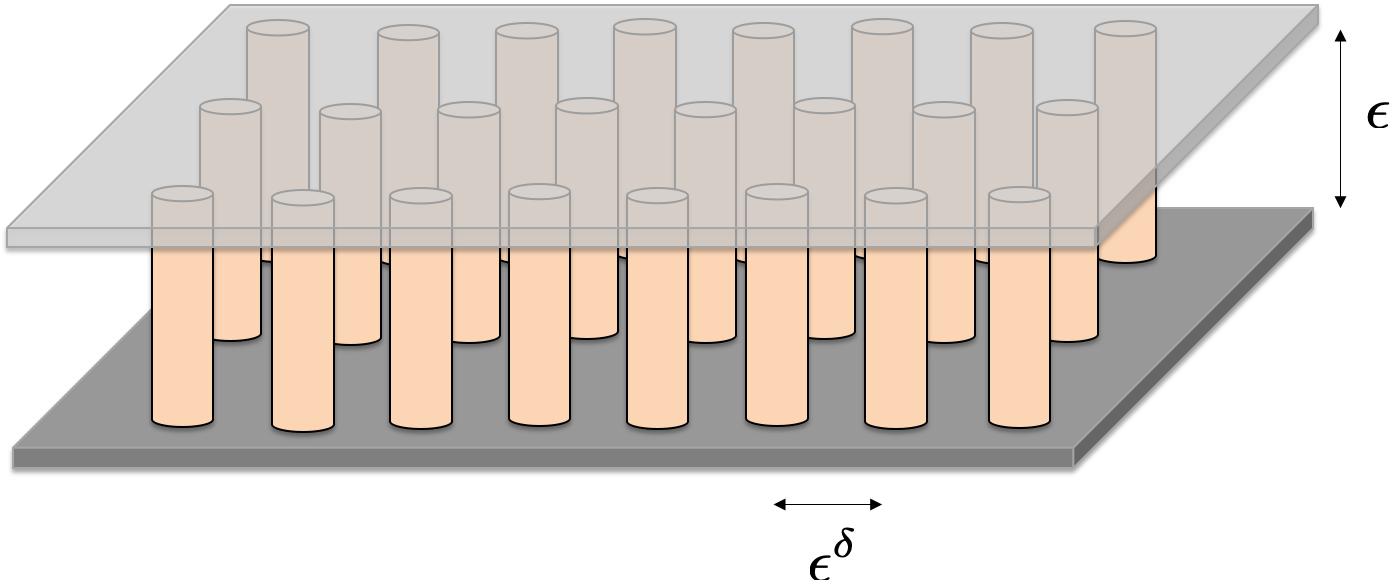}\hspace{1cm}  \includegraphics[width=5cm]{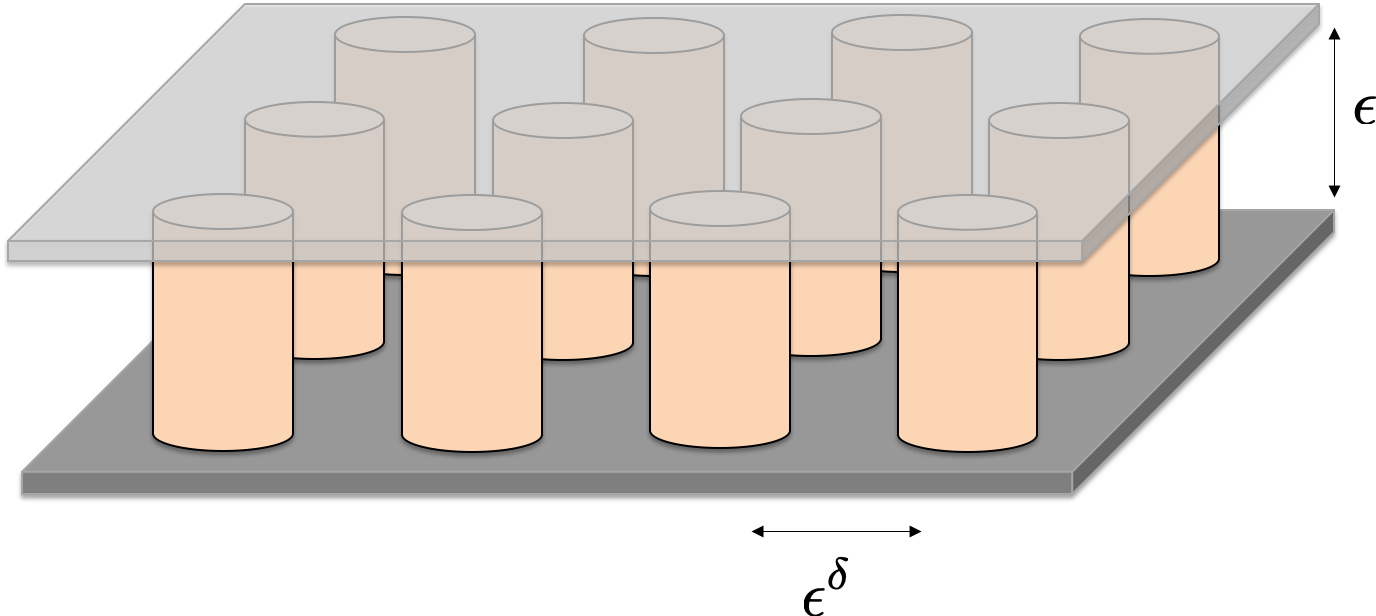}
\caption{View of {\bf HTPM} ($\delta>1$) (left) and {\bf PTPM} ($\delta=1$) (right)}
\label{HPTPM}
\end{figure}

\begin{figure}[h]
\centering \includegraphics[width=5cm]{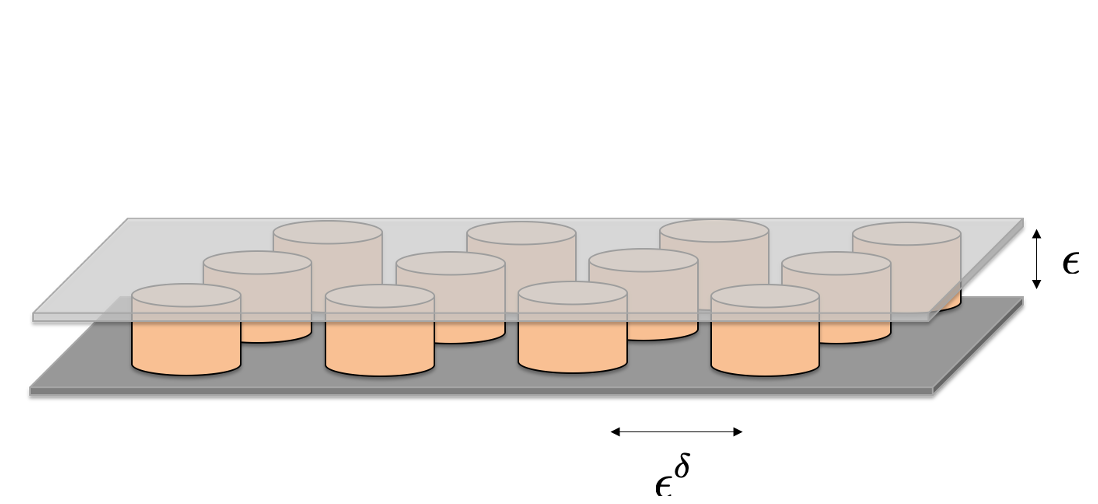}
\caption{View of {\bf VTPM} ($0<\delta<1$)}
\label{VTPM}
\end{figure}

Theoretical research concerning thin porous media ({\bf TPM}) is of great importance to application in various industries, see for instance Prat and Aga${\rm \ddot{e}}$sse \cite{Prat} and Yeghiazarian {\it et al.} \cite{Rosati}. The definition of thin porous media   is a domain whose thickness is small compared to the rest of the dimensions and which is heterogeneous, that is, it is perforated by periodically arranged cylindrical obstacles. If we denote the order of the thickness of the domain by the parameter $\epsilon$ and the order of the period and size of the diameter of the cylinders by $\epsilon^\delta$, with $\delta>0$, Fabricius {\it et al.} \cite{Fabricius} introduced  that for a Newtonian fluid flow through a thin porous media  that there is three different regimes depending on the relation of the thickness and the size of the cylinders:
\begin{itemize}
\item The homogeneously thin porous media {\bf (HTPM)}, corresponding to the case when the cylinder height is
much larger than interspatial distance, i.e. $\epsilon^\delta\ll \epsilon$, which is equivalent to $\delta>1$.  See Figure \ref{HPTPM} (left).

\item The proportionally thin porous media {\bf (PTPM)}, corresponding to the critical case when the cylinder height is proportional to the interspatial distance, i.e. $\epsilon^\delta\approx \epsilon$, which is equivalent to $\delta=1$. See Figure \ref{HPTPM} (right).

\item The very thin porous media {\bf (VTPM)}, corresponding to the case when the cylinder height is much
smaller than the interspatial distance, i.e. $\epsilon^\delta\gg \epsilon$,  which is equivalent to $0<\delta<1$. See Figure \ref{VTPM}.
\end{itemize}

\newpage 

Starting from the 3D Newtonian Stokes or Navier-Stokes system, it is deduced that the flow is governed by a 2D Darcy equation,  different in each regime, see Fabricius {\it et al.} \cite{Fabricius} for the proof by using the asymptotic expansion method and see Anguiano and Su\'arez-Grau \cite{Anguiano_SuarezGrau, Anguiano_SuarezGrau2} for a rigorous proof by using a version of the unfolding method.  We remark that in previous papers the Reynolds number is considered of order one and, since the thickness of the domain is small, the inertial term of the Navier-Stokes vanishes after the homogenization process  when $\epsilon$ tends to zero, which allows to derive the different 2D Darcy laws.\\

In  applications, it is important to know for what values of the Reynolds number the Darcy law models Newtonian fluid flows through a thin porous media. This issue is widely studied in porous media, where critical Reynolds numbers are obtained above which Darcy's laws are no longer valid, because the inertial term of the Navier-Stokes system must be taken into account, see for instance Dybbs and Edwards \cite{Dybbs}, Bourgeat and Mikeli\'c \cite{Bourgeat1}, Boettcher {\it et al.} \cite{Boettcher} or Zeng and Grigg \cite{ZZ1}.  In this sense, the only study related to this question in thin porous media is recently obtained by Jouybari and Lundstr$\ddot{\rm o}$m \cite{Jouybari}, where they present direct numerical simulations of post-Darcy flow in thin porous media. However, as far as the authors are aware, there are no theoretical studies where a critical Reynolds number is obtained below which Darcy's law is valid for each type of thin porous media described above, which is the goal and main novelty of this paper.\\

Thus, in this paper  we consider a Newtonian fluid flow governed by the Navier-Stokes system in a thin porous media $\Omega^\epsilon$. Assuming that the Reynolds number is proportional to the thickness of the domain and denoted by  $Re^\epsilon=\epsilon^{-\gamma}$, with $\gamma\in \mathbb{R}$, we consider the following system
\begin{equation}\label{N-S}
\left\{\begin{array}{rl}
-\eta\epsilon^\gamma  \Delta v^\epsilon + \nabla p^\epsilon +(v^\epsilon\cdot \nabla)v^\epsilon=f & \hbox{in }\Omega^\epsilon,\\
\\
{\rm div}\,v^\epsilon=0& \hbox{in }\Omega^\epsilon,
\\
\\
v^\epsilon=0& \hbox{on }\partial\Omega^\epsilon,
\end{array}\right.
\end{equation}
where $v^\epsilon$ is the velocity field, $p^\epsilon$ is the pressure, $\eta>0$ is the viscosity and $f$ is the external body  force.
\\

We prove that, for each regime of {\bf TPM}, there exists a critical value $\gamma_c$, i.e. a critical Reynolds number $Re_c=\epsilon^{-\gamma_c}$, such that for every $\gamma\leq \gamma_c$,  the fluid can be modeled by the corresponding Darcy law and above this quantity the inertial term has to be taken into account in the limit model. Observe that condition $\gamma\leq \gamma_c$ means that  the $Re$ is of order smaller or equal than  $\epsilon^{-\gamma_c}$. Thus, we obtain the following result:
\begin{itemize}
\item In the case {\bf HTPM} ($\delta>1$), the critical value is $\gamma_c=\delta$, i.e. the critical Reynolds number $Re_c=\epsilon^{-\delta}$.

\item In the case {\bf PTPM} ($\delta=1$), the critical value is $\gamma_c=1$, i.e. the critical Reynolds number  $Re_c=\epsilon^{-1}$.

 \item In the case {\bf VTPM} ($0<\delta<1$), the critical value is $\gamma_c=1$, i.e. the critical Reynolds number  $Re_c=\epsilon^{-1}$.
\end{itemize}

We remark that in Bourgeat and Mikeli\'c \cite{Bourgeat1}, they obtained a critical Reynolds number corresponding to $\gamma_c=3/2$ for Newtonian fluids through porous media with characteristic size $\epsilon$, and Boughanim and  Tapi\'ero \cite{Tapiero2} obtained the same critical value for Newtonian fluids through a thin slab of thickness $\epsilon$. In this work, we consider a thin heterogeneous media, which generalizes both studies. We summarize the order of the critical Reynolds numbers in the following table:
\begin{table}[h]\centering
\begin{tabular}{|c|c|c|c|c|}
\hline
 {\small Thin Slab}   &  {\small Porous media}  & {\small Homogeneously TPM} & {\small Proportionally TPM} &{\small Very TPM}  \\
   \hline\hline
  $ {\epsilon^{-{3\over 2}}}$ & $ {\epsilon^{-{3\over 2}}}$ & $ {\epsilon^{-\delta}}$ & ${\epsilon^{-1}}$ & $ {\epsilon^{-1}}$ \\
  \hline
\end{tabular} 
\caption{Critical Reynold number for Newtonian fluids depending on the type of {\bf TPM}.}
\label{table_Reynold}
\end{table}

  To prove this result, we first apply a dilatation in the vertical variable $x_3$, i.e. we introduce $z_3=x_3/\epsilon$, which let us to introduce the unknowns $(\tilde v^\epsilon,\tilde p^\epsilon)$ defined in a perforated domain with fixed thickness and denoted by $\widetilde\Omega^\epsilon$. Next, we need to prove that there exists an extension of the pressure, denoted by $\tilde P^\epsilon$, to the whole domain without obstacles $\Omega$. This will be done by duality arguments by means of a restriction operator introduced in Anguiano and Su\'arez-Grau \cite{Anguiano_SuarezGrau}.  Next, we have to prove that this extension is bounded in suitable Sobolev spaces (not only in $L^2$ as usual) depending on the {\bf TPM} considered.  It happens that estimating the extension of the pressure reduces to obtain sharp estimates for the inertial term of the dilated Navier-Stokes system depending on $\delta$ and $\gamma$. This will be made by using advanced interpolation techniques, which allow us to find the critical value $\gamma$ for each {\bf TPM}. We remark that using standard arguments, we would obtain that the extension of the pressure is bounded in the space $L^2$ and the critical Reynolds numbers would correspond to  $\gamma_c=3\delta/4$ for the {\bf HTMP} and $\gamma_c=3/4$ for the {\bf PTPM} and {\bf VTPM}. However, thanks to the sharp estimates obtained for the inertial term, we get a wider range of values of the Reynolds number for which the Darcy laws are valid, see Theorem \ref{mainthmEstimates} and its proof for more details. \\
  
\noindent  On the other hand, we complete the paper by obtaining the 2D Darcy laws for each {\bf TPM}. We obtain that the averaged velocity $\tilde V^\epsilon=\int_0^1\tilde v^\epsilon\,dz_3$ and the extension of the pressure $\tilde P^\epsilon$ satisfy the following approximations:
$$(\tilde V^\epsilon, \tilde P^\epsilon)\approx \left\{\begin{array}{l}  
  (\epsilon^{2\delta-\gamma} \tilde V, \tilde P) \quad  \hbox{for {\bf HTPM}},\\
\\
(\epsilon^{2-\gamma}\tilde V,  \tilde P)\quad  \hbox{for {\bf PTPM} and {\bf VTPM}},
 \end{array}\right. $$
where $(\tilde V, \tilde P)$, with $\tilde V_3=0$,  is the solution of a 2D Darcy law depending on the {\bf TPM} (see Theorems \ref{mainthmHTPM}, \ref{mainthmPTPM} and \ref{mainthmVTPM} for more details):
$$
\begin{array}{c}
\displaystyle \tilde V'(x')={1\over \eta}\mathcal{K}\left(f'(x')-\nabla_{x'}\tilde P(x')\right),\quad \tilde V_3(x')=0\quad\hbox{in}\ \omega,
\medskip\\
\displaystyle
{\rm div}_{x'}\tilde V'(x')=0\quad\hbox{in}\ \omega,\quad \tilde V'(x')\cdot n=0\quad\hbox{on}\ \partial\omega,
\end{array}
$$
where $n$ is the outer normal vector to $\partial\omega$. We remark that in the {\bf HTPM} the permeability tensor $\mathcal{K}\in \mathbb{R}^{2\times 2}$ is calculated by solving 2D Stokes local problems, in the {\bf PTPM}  by solving 3D Stokes local problems, and in the {\bf VTPM}  by solving 2D Hele-Shaw local problems. We also remark that the approximations depend on  the Reynolds parameter $\gamma$ and the type of {\bf TPM}. The proof of this result follows the proof of Anguiano and Su\'arez-Grau \cite{Anguiano_SuarezGrau, Anguiano_SuarezGrau2} by taking  the Reynolds number $\epsilon^{-\gamma}$ into account in the estimates of the velocity. Since the novelty of this paper is the introduction of the Reynolds number and the inertial term in the Navier-Stokes system, in the homogenization process when $\epsilon$ tends to zero we focus our attention in proving that the inertial term vanishes in the limit by using the sharp estimates of the pressure obtained previously for each type of {\bf TPM}. In order to have a self contained study, we also give details of how obtaining the corresponding 2D Darcy laws. \\

As far as the authors know, this is the first attempt to carry out such a theoretical analysis for Newtonian fluids in {\bf TPM}, which could be instrumental for understanding the effects of the Reynolds number on the behavior of Newtonian fluid flows in each type of {\bf TPM}. In view of that, more efficient numerical algorithms could be developed improving, hopefully, the known engineering practice.\\

To finish the introduction, we give a list of recent references concerning studies of partial differential equations in {\bf TPM}.  Studies related to Newtonian fluids through {\bf TPM} can be found in Anguiano and Su\'arez-Grau \cite{Anguiano_SuarezGrau_Derivation, Anguiano_SG_Net}, Armiti-Juber \cite{Armiti}, Bayada and {\it et al.}  \cite{Bayada}, Larsson {\it et al.}  \cite{Sofia}, Su\'arez-Grau \cite{SG_MN}, Valizadeh and   Rudman \cite{Valizadeh},   Wagner {\it et al.} \cite{Wagner} and Zhengan and  Hongxing \cite{ZZ}. Concerning generalized Newtonian fluids see Anguiano and Su\'arez-Grau \cite{Carreau_Ang_Bonn_SG, Anguiano_SuarezGrau_CMS, AnguianoSG_small}, for Bingham fluids see Anguiano and Bunoiu \cite{Ang-Bun, Ang-Bun2}, for compressible and piezo-viscous flow  see P\'erez-R\`afols {\it et al.}  \cite{Perez}, and for micropolar fluids  see Su\'arez-Grau \cite{SG_micropolar} and  for diffusion problems  see Anguiano \cite{Anguiano_BMMS, Anguiano_plaplace} and Bunoiu and Timofte \cite{Bunoiu}.\\

The paper is organized as follows. In Section \ref{Sec:2}, we introduce the description of the {\bf TPM} considered in this paper. In Section \ref{sec:Estimates},  we introduce the problem under consideration and we obtain the first main result concerning the sharp estimates for the extension of the pressure, see Theorem \ref{mainthmEstimates}. In Section \ref{sec:Darcy}, we give describe the Darcy laws for each type of {\bf TPM}, see Theorems \ref{mainthmHTPM}, \ref{mainthmPTPM} and \ref{mainthmVTPM}. We give the proof of previous results  in Section \ref{sec:Proofs} and we finish the paper with a list of references.

 \section {Description of the TPM}\label{Sec:2} 
 
Let us introduce some notation that will be useful in the rest of the paper. The points $x\in\mathbb{R}^3$ will be decomposed as $x=(x^{\prime},x_3)$ with $x^{\prime}=(x_1,x_2)\in \mathbb{R}^2$, $x_3\in \mathbb{R}$. We also use the notation $x^{\prime}$ to denote a generic vector of $\mathbb{R}^2$. 

\paragraph{Description of the {\bf TPM}.} The periodic porous medium is defined by a domain $\omega$ and an associated microstructure, or periodic cell $Z^{\prime}=(-1/2,1/2)^2$, which is made of two complementary parts: the fluid part $Z^{\prime}_{f}$, and the solid part $S^{\prime}$ ($Z^{\prime}_f  \bigcup S^{\prime}=Z^\prime$ and $Z^{\prime}_f  \bigcap S^{\prime}=\emptyset$). More precisely, we assume that $\omega$ is a smooth, bounded, connected set in $\mathbb{R}^2$, and that $S^{\prime}$ is an open connected subset of $Z^\prime$ with a smooth boundary $\partial S^\prime$, such that $\bar S^\prime$ is strictly included  in $Z^\prime$.\\

The microscale of a porous medium is denoted by ${\epsilon}^\delta$, with $\delta\in \mathbb{R}$. The domain $\omega$ is covered by a regular mesh of square of size ${\epsilon^\delta}$: for $k^{\prime}\in \mathbb{Z}^2$, each cell $Z^{\prime}_{k^{\prime},{\epsilon}^\delta }={\epsilon}^\delta k^{\prime}+{\epsilon}^\delta Z^{\prime}$ is divided in a fluid part $Z^{\prime}_{f_{k^{\prime}},{\epsilon}^\delta }$ and a solid part $S^{\prime}_{k^{\prime},{\epsilon}^\delta }$, i.e. is similar to the unit cell $Z^{\prime}$ rescaled to size ${\epsilon}^\delta$. We define $Z=Z^{\prime}\times (0,1)\subset \mathbb{R}^3$, which is divided in a fluid part $Z_{f}=Z'_f\times (0,1)$ and a solid part $S=S'\times(0,1)$, and consequently $Z_{k^{\prime},{\epsilon}^\delta }=Z^{\prime}_{k^{\prime},{\epsilon}^\delta }\times (0,1)\subset \mathbb{R}^3$, which is also divided in a fluid part $Z_{f_{k^{\prime}},{\epsilon}^\delta }$ and a solid part $S_{{k^{\prime}},{\epsilon}^\delta }$.\\

We denote by $\tau(\overline S'_{k',\epsilon^\delta} )$ the set of all translated images of $\overline S'_{k',\epsilon^\delta }$. The set $\tau(\overline S'_{k',\epsilon^\delta })$ represents the obstacles in $\mathbb{R}^2$.\\

The fluid part of the bottom $\omega^{\epsilon}\subset \mathbb{R}^2$ of a porous medium is defined by $\omega^{\epsilon}=\omega\backslash\bigcup_{k^{\prime}\in \mathcal{K}_{\epsilon^\delta} } \bar S^{\prime}_{{k^{\prime}},{\epsilon^\delta} },$ where $\mathcal{K}^{\epsilon}=\{k^{\prime}\in \mathbb{Z}^2: Z^{\prime}_{k^{\prime}, {\epsilon}^\delta } \cap \omega \neq \emptyset \}$.  The whole fluid part $\Omega^{\epsilon}\subset \mathbb{R}^3$ in the thin porous medium is defined by 
\begin{equation}\label{Dominio1}
\Omega^{\epsilon}=\{  (x', x_3)\in \omega^{\epsilon}\times \mathbb{R}: 0<x_3<\epsilon \}.
\end{equation}
We assume that the obstacles $\tau(\overline S'_{k',\epsilon^\delta })$ do not intersect the boundary $\partial\omega$ and we denote by $\mathbb{S}^\epsilon$ the set of the solid cylinders contained in $\Omega^\epsilon$, i.e. $\mathbb{S}^\epsilon=\bigcup_{k^{\prime}\in \mathcal{K}_{\epsilon}^\delta } S'_{k^\prime, \epsilon^\delta }\times (0,\epsilon)$.
\\

\noindent We define  the following sets:
\begin{itemize}
\item Dilated TPM:
\begin{equation}\label{OmegaTilde}
\widetilde{\Omega}^{\epsilon}=\omega^{\epsilon}\times (0,1),
\end{equation}
\item Defining $\mathbb{T}^\epsilon=\bigcup_{k^{\prime}\in \mathcal{K}_{\epsilon}} T_{k^\prime, \epsilon^\delta }$ as the set of the solid cylinders contained in $\widetilde \Omega^\epsilon$, it holds
$$\widetilde{\Omega}_{\epsilon}=\Omega\backslash \mathbb{T}^\epsilon.$$

\item Domains without perforations:
\begin{equation}\label{OmQep}
\Omega=\omega\times (0,1), \quad Q^\epsilon=\omega\times (0,\epsilon).
\end{equation}
 \end{itemize}

 \section {Sharp pressure estimates in {\bf TPM}}\label{sec:Estimates}
 In this section, we give the main result concerning the sharp pressure estimates in {\bf TPM} and the critical Reynolds number in each type of {\bf TPM}.  To do this, we consider a Newtonian fluid flow governed by the Navier-Stokes system (\ref{N-S}) in the thin porous media $\Omega^\epsilon$.  Assuming that the second member $f$ of (\ref{N-S}) is of the form
\begin{equation}\label{fassump}
f (x)=(f'(x'),0)\quad \hbox{with }f'\in L^2(\omega)^2,
\end{equation}
it is known that there exists at least one weak solution $(v^\epsilon, p^\epsilon)\in H^1_0(\Omega^\epsilon)\times L^2_0(\Omega^\epsilon)$ of (\ref{N-S}), where $L_0^2$ is the space of functions of $L^2$  with null mean value.\\

\noindent In order to work in a domain with fixed height,  we use the dilatation in the variable $x_3$ as follows
\begin{equation}\label{dilatacion}
z_3=\frac{x_3}{\epsilon},
\end{equation}
to have the functions defined in $\widetilde\Omega^{\epsilon}$. This leads to define $\tilde{v}^{\epsilon}$ and $\tilde{p}^{\epsilon }$ by $$\tilde{v}^{\epsilon }(x^{\prime},z_3)=v^{\epsilon }(x^{\prime}, \epsilon z_3),\text{\ \ }\tilde{p}^{\epsilon }(x^{\prime},z_3)=p^{\epsilon }(x^{\prime}, \epsilon z_3), \text{\ \ } a.e.\text{\ } (x^{\prime},z_3)\in \widetilde{\Omega}^{\epsilon }.$$
Then, applying this change of variables to the system (\ref{N-S}), we consider the following dilated Navier-Stokes system
\begin{equation}\label{N-S-d}
\left\{\begin{array}{rl}
-\eta\epsilon^\gamma\Delta_\epsilon \tilde v^\epsilon + \nabla_\epsilon\, \tilde p^\epsilon +(\tilde v^\epsilon\cdot \nabla_\epsilon)\tilde v^\epsilon=f & \hbox{in }\widetilde \Omega^\epsilon,\\
\\
{\rm div}_\epsilon\,\tilde v^\epsilon=0& \hbox{in }\widetilde \Omega^\epsilon,
\\
\\
\tilde v^\epsilon=0& \hbox{on }\partial\widetilde \Omega^\epsilon,
\end{array}\right.
\end{equation}
where $\Delta_\epsilon\tilde v^\epsilon=\Delta_{x'} \tilde v^\epsilon+\epsilon^{-2}\partial_{z_3}^2 \tilde v^\epsilon$, $\nabla_\epsilon\tilde v^\epsilon=(\nabla_{x'}\tilde v^\epsilon, \epsilon^{-1}\partial_{z_3}\tilde v^\epsilon)^t$, and  ${\rm div}_\epsilon\,\tilde v^\epsilon={\rm div}_{x'}(\tilde v^\epsilon)'+\epsilon^{-1}\tilde v^{\epsilon}_3$.
\\

Below, we give the main result concerning the existence of an extension of the pressure $\tilde p^\epsilon$ and the sharp estimates of the pressure depending on $\delta$ and $\gamma$, which will give the existence of critical Reynolds number in each type of {\bf TPM}.  To prove this, we take into account the estimates of the velocity and we prove that  there exists an extension $\tilde P^\epsilon$ of the pressure $\tilde p^\epsilon$ bounded in a Sobolev space $L^{C(\delta)}$   (not only in $L^2$ as usual) for $C(\delta)\in (1,2]$ being the conjugate of a certain $r$, $r\geq 2$. The existence of an extension of the pressure is proved by means of duality argument  by using a restriction operator $\tilde R^\epsilon_r$ acting from $W^{1,r}_0(\Omega)^3$ into $W^{1,r}_0(\widetilde \Omega^\epsilon)^3$,  which is introduced in \cite{Anguiano_SuarezGrau, Anguiano_SuarezGrau2}. In order  to prove that the extension of the pressure is bounded and that inertial term of the dilated Navier-Stokes system (\ref{N-S-d}) vanishes in the limit, we have to prove that
$$\left|\int_{\widetilde\Omega^\epsilon}(\tilde v^\epsilon\cdot \nabla_\epsilon)v^\epsilon\tilde R^\epsilon_r(\tilde \varphi)\,dx'dz_3\right|\leq C\epsilon^\alpha\|\tilde \varphi\|_{W^{1,r}_0(\Omega)^3},\quad \forall\,\tilde\varphi\in W^{1,r}_0(\Omega)^3,$$
with $\alpha\geq 0$. This inequality will provide the sharp estimates of the extension of the pressure and the range of values of $\gamma$ for which the inertial term will vanish in the limit. \\

\begin{theorem}[Sharp pressure estimates]\label{mainthmEstimates}
Consider  any weak solution $\tilde p^\epsilon$ of (\ref{N-S-d}) and $C(\delta)$ defined by
\begin{equation}\label{C_case_H}
C(\delta)=\left\{\begin{array}{ccl}
2&\hbox{if}& \gamma\leq \displaystyle{3\delta\over 4},
\\
\\
\displaystyle {3\delta\over 2\gamma}&\hbox{if}& \displaystyle{3\delta\over 4}<\gamma\leq  \delta.
\end{array}\right.
\end{equation}
\\

\noindent Then,   depending on the type of {\bf TPM}:
\begin{itemize}
\item[(i)] In the case {\bf HTPM} ($\delta>1$) and $\gamma \leq \delta$,  there exist an extension $\tilde P^\epsilon\in L^{C(\delta)}_0( \Omega)$ of the pressure $\tilde p^\epsilon$ satisfying (\ref{N-S-d}) and a positive constant $C$, independent of $\epsilon$, such that 
\begin{equation}\label{esti_p_caseH}
\|\tilde P^\epsilon\|_{L^{C(\delta)}(\Omega)}\leq C,\quad \|\nabla_\epsilon \tilde P^\epsilon\|_{W^{-1, C(\delta)}(\Omega)^3}\leq C.
\end{equation}
\item[(ii)] In the case {\bf PTPM} ($\delta=1$) and $\gamma\leq 1$, there exist an extension $\tilde P^\epsilon\in L^{C(1)}_0( \Omega)$ of  the pressure $\tilde p^\epsilon$ and a positive constant $C$, independent of $\epsilon$, such that 
\begin{equation}\label{esti_p_caseP}
\|\tilde P^\epsilon\|_{L^{C(1)}(\Omega)}\leq C,\quad \|\nabla_\epsilon \tilde P^\epsilon\|_{W^{-1, C(1)}(\Omega)^3}\leq C.
\end{equation}
\item[(iii)] In the case {\bf VTPM} ($0<\delta<1$) and $\gamma\leq 1$, there exist an extension $\tilde P^\epsilon\in L^{C(1)}_0( \Omega)$ of  the pressure $\tilde p^\epsilon$ and a positive constant $C$, independent of $\epsilon$, such that 
\begin{equation}\label{esti_p_caseV}
\|\tilde P^\epsilon\|_{L^{C(1)}(\Omega)}\leq C,\quad \|\nabla_\epsilon \tilde P^\epsilon\|_{W^{-1, C(1)}(\Omega)^3}\leq C.
\end{equation}\end{itemize}
\end{theorem}

\section{Darcy's laws for each type of  thin porous media}\label{sec:Darcy}
In this section, we take advantage of the results of Theorem \ref{mainthmEstimates} to show that the inertial term disappears in the limit and we obtain the lower-dimensional Darcy laws for each type of {\bf TPM}. This will be done by means of an adaptation (see  \cite{Anguiano_SuarezGrau, Anguiano_SuarezGrau2}) of the unfolding method \cite{CDG}.
\\ 

 \begin{theorem}[Darcy's law for {\bf HTPM}]\label{mainthmHTPM}
Consider $\delta>1$,   $\gamma \leq \delta$ and $C(\delta)$ defined by (\ref{C_case_H}).  Then, there exists $\tilde v\in L^2(\Omega)^3$ with $\tilde v_3=0$ and $\tilde P\in L_0^{C(\delta)}(\omega)$, such that the extension $(\tilde v^\epsilon, \tilde P^\epsilon)$ of a solution of (\ref{N-S-d}) satisfies the convergences
$$\epsilon^{\gamma-2\delta}\tilde v^\epsilon\rightharpoonup \tilde v\quad\hbox{in}\ L^2(\Omega)^3,\quad \tilde P^\epsilon\to \tilde P\quad\hbox{in}\ L^{C(\delta)}(\Omega).
$$
Moreover, defining $\tilde V(x')=\int_0^1\tilde v(x',z_3)\,dz_3$, we have that $(\tilde V, \tilde P)\in L^2(\omega)^3\times (L^2_0(\omega)\times H^1(\omega))$ is the unique solution of the linear Darcy law 
\begin{equation}\label{DarcyHTPM}
\begin{array}{c}
\displaystyle \tilde V'(x')={1\over \eta}\mathcal{K}_H\left(f'(x')-\nabla_{x'}\tilde P(x')\right),\quad \tilde V_3(x')=0\quad\hbox{in}\ \omega,
\medskip\\
\displaystyle
{\rm div}_{x'}\tilde V'(x')=0\quad\hbox{in}\ \omega,\quad \tilde V'(x')\cdot n=0\quad\hbox{on}\ \partial\omega.
\end{array}
\end{equation}
Here, $\mathcal{K}_H\in \mathbb{R}^{2\times 2}$ is the permeability tensor, which is symmetric
and positive definite,  defined by 
\begin{equation}\label{JHTPM}
(\mathcal{K}_H)_{ij}=\int_{Z'_f}w^i_j(z')\,dz',\quad i,j=1,2,
\end{equation}
with $w^i \in H^1_{\#}(Z'_f)^2$, $i=1,2$, the unique solution of the local 2D Stokes problem defined by
\begin{equation}\label{cell_problem_HTPM}
\left\{\begin{array}{rl}
-\Delta_{z'}w^i + \nabla_{z'}\pi^i=e_i &\hbox{in}\  Z_f',
\medskip\\
{\rm div}_{z'} w^i=0&\hbox{in}\  Z_f',
\medskip\\
w^i=0&\hbox{in}\  S',
\medskip\\
w^i, \pi^i\ Z'-\hbox{periodic}.
\end{array}
\right.
\end{equation}\\
 \end{theorem}

 \begin{theorem}[Darcy's law for {\bf PTPM}]\label{mainthmPTPM}
Consider $\delta=1$, $\gamma\leq 1$ and $C(\delta)$ defined by (\ref{C_case_H}).  Then, there exists $\tilde v\in H^1(0,1;L^2(\omega)^3)$ with $\tilde v_3=0$ and $\tilde v=0$ on $z_3=\{0,1\}$, and $\tilde P\in L_0^{C(1)}(\omega)$, such that the extension $(\tilde v^\epsilon, \tilde P^\epsilon)$ of a solution of (\ref{N-S-d}) satisfies the convergences
$$\epsilon^{\gamma-2}\tilde v^\epsilon\rightharpoonup \tilde v\quad\hbox{in}\ H^1(0,1;L^2(\omega)^3),\quad \tilde P^\epsilon\to \tilde P\quad\hbox{in}\ L^{C(1)}(\omega).
$$
Moreover, defining $\tilde V(x')=\int_0^1\tilde v(x',z_3)\,dz_3$, we have that $(\tilde V, \tilde P)\in L^2(\omega)^3\times (L^2_0(\omega)\times H^1(\omega))$ is the unique solution of the linear Darcy law 
\begin{equation}\label{DarcyPTPM}
\begin{array}{c}
\displaystyle \tilde V'(x')={1\over \eta}\mathcal{K}_P\left(f'(x')-\nabla_{x'}\tilde P(x')\right),\quad \tilde V_3(x')=0\quad\hbox{in}\ \omega,
\medskip\\
\displaystyle
{\rm div}_{x'}\tilde V'(x')=0\quad\hbox{in}\ \omega,\quad \tilde V'(x')\cdot n=0\quad\hbox{on}\ \partial\omega.
\end{array}
\end{equation}
Here, $\mathcal{K}_P\in \mathbb{R}^{2\times 2}$ is the permeability tensor, which is symmetric
and positive definite,  defined by 
\begin{equation}\label{JPTPM}
(\mathcal{K}_P)_{ij}=\int_{Z_f}w^i_j(z)\,dz,\quad i,j=1,2,
\end{equation}
with $w^i \in H^1_{\#}(Z_f)^3$, $i=1,2$, the unique solution of the local 3D Stokes problem defined by
\begin{equation}\label{cell_problem_PTPM}
\left\{\begin{array}{rl}
-\Delta_{z}w^i + \nabla_{z}\pi^i=e_i &\hbox{in}\  Z_f,
\medskip\\
{\rm div}_{z} w^i=0&\hbox{in}\  Z_f,
\medskip\\
w^i=0&\hbox{in}\  S,
\medskip\\
w^i=0&\hbox{on}\  z_3=\{0,1\}.
\end{array}
\right.
\end{equation}\\
 \end{theorem}

  \begin{theorem}[Darcy's law for {\bf VTPM}]\label{mainthmVTPM}
Consider $\delta\in (0,1)$, $\gamma\leq 1$ and $C(\delta)$ defined by (\ref{C_case_H}).  Then, there exists $\tilde v\in H^1(0,1;L^2(\omega)^3)$ with $\tilde v_3=0$ and $\tilde v=0$ on $z_3=\{0,1\}$, and $\tilde P\in L_0^{C(1)}(\omega)$, such that the extension $(\tilde v^\epsilon, \tilde P^\epsilon)$ of a solution of (\ref{N-S-d}) satisfies the convergences
$$\epsilon^{\gamma-2}\tilde v^\epsilon\rightharpoonup \tilde v\quad\hbox{in}\ H^1(0,1;L^2(\omega)^3),\quad \tilde P^\epsilon\to \tilde P\quad\hbox{in}\ L^{C(1)}(\omega).
$$
Moreover, defining $\tilde V(x')=\int_0^1\tilde v(x',z_3)\,dz_3$, we have that $(\tilde V, \tilde P)\in L^2(\omega)^3\times (L^2_0(\omega)\times H^1(\omega))$ is the unique solution of the linear Darcy law 
\begin{equation}\label{DarcyVTPM}
\begin{array}{c}
\displaystyle \tilde V'(x')={1\over 12\eta}\mathcal{K}_V\left(f'(x')-\nabla_{x'}\tilde P(x')\right),\quad \tilde V_3(x')=0\quad\hbox{in}\ \omega,
\medskip\\
\displaystyle
{\rm div}_{x'}\tilde V'(x')=0\quad\hbox{in}\ \omega,\quad \tilde V'(x')\cdot n=0\quad\hbox{on}\ \partial\omega.
\end{array}
\end{equation}
Here, $\mathcal{K}_V\in \mathbb{R}^{2\times 2}$ is the permeability tensor, which is symmetric
and positive definite,  defined by 
\begin{equation}\label{JVTPM}
(\mathcal{K}_V)_{ij}=\int_{Z'_f}\left(\nabla_{z'}\pi^i(z')+e_i\right)e_j\,dz',\quad i,j=1,2,
\end{equation}
with $\pi^i\in  H^1_{\#}(Z'_f)$, $i=1,2$, the unique solution of the local Hele-Shaw problem defined by
\begin{equation}\label{cell_problem_VTPM}
\left\{\begin{array}{rl}
-\Delta_{z'}\pi^i=0 &\hbox{in}\  Z'_f,
\medskip\\
(\nabla_{z'}\pi^i+e_i)\cdot n=0&\hbox{on}\  \partial S'.
\end{array}
\right.
\end{equation}
 \end{theorem}
\section{Proof of the results}\label{sec:Proofs}
In this section we provide the proof of Theorems \ref{mainthmEstimates}, \ref{mainthmHTPM}, \ref{mainthmPTPM} and \ref{mainthmVTPM}. In Subsection \ref{subsec:VelEstim}, we give the estimates of the velocity and the extension of the velocity. In Subsection \ref{Subsec:Pressure}, we introduce the restriction operator to extend the pressure to the whole domain and we give the proof of Theorem \ref{mainthmEstimates}. We introduce the adaptation of the unfolding method in Subsection \ref{Subsec:unfolded} and we give the variational formulation in Subsection \ref{Subsec:varfun}. The proof of Theorem \ref{mainthmHTPM} is given in Subsection \ref{Subsec:HTPM}, the proof of Theorem \ref{mainthmPTPM} in Subsection \ref{Subsec:PTPM}, and the proof of Theorem \ref{mainthmVTPM} in Subsection \ref{Subsec:VTPM}.
\subsection{Velocity estimates}\label{subsec:VelEstim}
First, we deduce the velocity estimates. To do this, let us first recall  the  Poincar\'e inequality  in $\widetilde\Omega^\epsilon$, see \cite[Lemma 4.2 and Remark 4.3]{Anguiano_SuarezGrau} for more details.
\\

\begin{lemma}[Poincar\'e's inequality] There exists a positive constant $C$, independent of $\epsilon$, such that
\begin{itemize}
\item   In the case {\bf HTPM}  ($\delta>1$), it holds
\begin{equation}\label{PoincareH}
\|\tilde \varphi\|_{L^2(\widetilde \Omega^\epsilon)^3}\leq C\epsilon^\delta \|D_\epsilon \tilde  \varphi\|_{L^2(\widetilde \Omega^\epsilon)^{3\times 3}},\quad \forall\,\tilde  \varphi\in H^1_0(\widetilde \Omega^\epsilon)^3.
\end{equation}
\item In the case {\bf PTPM}  ($\delta=1$), it holds
\begin{equation}\label{PoincareP}
\|\tilde \varphi\|_{L^2(\widetilde \Omega^\epsilon)^3}\leq C\epsilon  \|D_\epsilon \tilde  \varphi\|_{L^2(\widetilde \Omega^\epsilon)^{3\times 3}},\quad \forall\,\tilde  \varphi\in H^1_0(\widetilde \Omega^\epsilon)^3.
\end{equation}
\item In the case  {\bf VTPM}  ($0<\delta<1$), it holds
\begin{equation}\label{PoincareV}
\|\tilde \varphi\|_{L^2(\widetilde \Omega^\epsilon)^3}\leq C\epsilon \|D_\epsilon \tilde  \varphi\|_{L^2(\widetilde \Omega^\epsilon)^{3\times 3}},\quad \forall\,\tilde  \varphi\in H^1_0(\widetilde \Omega^\epsilon)^3.
\end{equation}
\end{itemize}
\end{lemma}
By using the previous lemma and taking into account that  it appears the Reynolds number $\epsilon^\gamma$ in system (\ref{N-S-d}), the following velocity estimates depending on the type of the thin porous media are straightforward obtained by following the steps of the proof of \cite[Lemma 4.4]{Anguiano_SuarezGrau}. 
\\

\begin{lemma}[Velocity estimates]\label{Vel_Estimates_lemma} Let $\tilde v^\epsilon$ a weak  solution of (\ref{N-S-d}). Then:
\begin{itemize}
\item[(i)] In the case  {\bf HTPM} ($\delta>1$),  there exists a positive constant $C$, independent of $\epsilon$, such that we have
\begin{equation}\label{estimates_v_H}
\|\tilde v^\epsilon\|_{L^2(\widetilde \Omega^\epsilon)^{3}}\leq C\epsilon^{2\delta-\gamma} ,\quad \|D_\epsilon \tilde v^\epsilon\|_{L^2(\widetilde \Omega^\epsilon)^{3\times 3}}\leq C\epsilon^{\delta-\gamma} \,.
\end{equation}
\item[(ii)] In the case {\bf PTPM} ($\delta=1$),  there exists a positive constant $C$, independent of $\epsilon$, such that we have
\begin{equation}\label{estimates_v_P}
\|\tilde v^\epsilon\|_{L^2(\widetilde \Omega^\epsilon)^{3}}\leq C\epsilon^{2-\gamma} ,\quad \|D_\epsilon \tilde v^\epsilon\|_{L^2(\widetilde \Omega^\epsilon)^{3\times 3}}\leq C\epsilon^{1-\gamma} \,.
\end{equation}
\item[(iii)] In the case {\bf VTPM} ($0<\delta<1$),  there exists a positive constant $C$, independent of $\epsilon$, such that we have
\begin{equation}\label{estimates_v_V}
\|\tilde v^\epsilon\|_{L^2(\widetilde \Omega^\epsilon)^{3}}\leq C\epsilon^{2-\gamma} ,\quad \|D_\epsilon \tilde v^\epsilon\|_{L^2(\widetilde \Omega^\epsilon)^{3\times 3}}\leq C\epsilon^{1-\gamma} \,.
\end{equation}
\end{itemize}
\end{lemma}
\begin{remark}[Extension of velocity to $\Omega$] We make the following remarks:
\begin{itemize}
\item[--] The velocity $\tilde v^\epsilon$ will be extended  by zero in $\Omega\setminus\widetilde \Omega^\epsilon$, which  is compatible with the zero boundary condition on the boundary of the obstacles and the exterior boundary of $\Omega$. 
\item[--] We will denote the extension of the velocity by same symbol. 

\item[--] The extended velocity is divergence free.

\item[--] The estimates given in Lemma \ref{Vel_Estimates_lemma} are still valid for the extended velocity.
\end{itemize}
\end{remark}
\subsection{Extension of pressure and proof of Theorem \ref{mainthmEstimates}}\label{Subsec:Pressure}
\noindent To extend the pressure to the whole domain $\Omega$, we use the mapping $R^\epsilon_q$, $1<q<+\infty$, defined in \cite[Lemma 4.5]{Anguiano_SuarezGrau}, which allows us to extend the pressure $p^\epsilon$ from $\Omega^\epsilon$ to $Q^\epsilon$. 
\\

\begin{lemma}[Lemma 4.5-(i) in \cite{Anguiano_SuarezGrau}] \label{restriction_operator}
There exists  a (restriction) operator $R^\epsilon_r$ acting from $W^{1,r}_0(Q^\epsilon)^3$ into $W^{1,r}_0(\Omega^\epsilon)^3$, $1<r<+\infty$, such that
\begin{enumerate}
\item $R^\epsilon_r( \varphi)=\varphi$, if $\varphi \in W^{1,r}_0(\Omega^\epsilon)^3$ (elements of $W^{1,r}_0(\Omega^\epsilon)$ are extended by $0$ to $Q^\epsilon$).
\item ${\rm div}R^\epsilon_r( \varphi)=0\hbox{  in }\Omega^\epsilon$, if ${\rm div}\,\varphi=0\hbox{  on }Q^\epsilon$.
\item For every $\varphi\in W^{1,r}_0(Q^\epsilon)^3$, there exists a positive constant $C$, independent of $\varphi$ and $\epsilon$, such that
\begin{itemize}
\item In the case {\bf HTPH} ($\delta>1$), it holds
\begin{equation}\label{estim_restricted_H}
\begin{array}{l}
\!\!\!\!\!\!\!\!\!\!\!\!\!\|R^\epsilon_r (\varphi)\|_{L^r(\Omega^\epsilon)^{3}}+ \epsilon^\delta \|D R^\epsilon_r( \varphi)\|_{L^r(\Omega^\epsilon)^{3\times 3}} \leq C\left(\|\varphi\|_{L^r(Q^\epsilon)^3}+\epsilon^\delta \|D \varphi\|_{L^r(Q^\epsilon)^{3\times 3}}\right).
\end{array}
\end{equation}
\item In the case {\bf PTPM} ($\delta=1$), it holds
\begin{equation}\label{estim_restricted_P}
\begin{array}{l}
\!\!\!\!\!\!\!\!\!\!\!\!\!\|R^\epsilon_r( \varphi)\|_{L^r(\Omega^\epsilon)^{3}}+ \epsilon\|D R^\epsilon_r( \varphi)\|_{L^r(\Omega^\epsilon)^{3\times 3}} \leq C\left(\|\varphi\|_{L^r(Q^\epsilon)^3}+\epsilon \|D \varphi\|_{L^r(Q^\epsilon)^{3\times 3}}\right).
\end{array}
\end{equation}
\item In the case {\bf VTPM} ($0<\delta<1$), it holds
\begin{equation}\label{estim_restricted_V}
\begin{array}{l}
\!\!\!\!\!\!\!\!\!\!\|R^\epsilon_r( \varphi)\|_{L^r(\Omega^\epsilon)^{3}}+ \epsilon\|D R^\epsilon_r( \varphi)\|_{L^r(\Omega^\epsilon)^{3\times 3}} \leq C\left(\|\varphi\|_{L^r(Q^\epsilon)^3}+\epsilon^\delta \|D \varphi\|_{L^r(Q^\epsilon)^{3\times 3}}\right).
\end{array}
\end{equation}
\end{itemize}
\end{enumerate}
\end{lemma}

\noindent Now, we are able to prove the main result concerning the extension of the pressure and the derivation of sharp estimates.\\

\noindent {\bf Proof of Theorem \ref{mainthmEstimates}}. 
We divide the proof in four steps. In the first step, we extend the pressure by duality arguments and in the rest of steps,  we obtain the estimates of the pressure in the case of {\bf HTPM}, {\bf PTPM} and {\bf VTPM}, respectively.\\

\noindent {\it Step 1.} {\it Extension of the pressure}.  In this step, we use the restriction operator  given in  Lemma \ref{restriction_operator}, i.e.  $R^\epsilon_r$, with $r\geq 2$ where $r$ will be determined in the following steps (depending on the case of {\bf TPM}).
\\

\noindent We  introduce $F^\epsilon$ in $W^{-1,C(\delta)}(Q^\epsilon)^3$, where $C(\delta)$ is the conjugate of $r$,  in the following way
\begin{equation}\label{F}\langle F^\epsilon, \varphi\rangle_{W^{-1,C(\delta)}(Q^\epsilon)^3, W^{1,r}_0(Q^\epsilon)^3}=\langle \nabla p^\epsilon, R^\epsilon_r(\varphi)\rangle_{{W^{-1,C(\delta)}(\Omega^\epsilon)^3, W^{1,r}_0(\Omega^\epsilon)^3}}\,,
\end{equation}
for every $\varphi\in W^{1,r}_0(Q^\epsilon)^3$. By using the variational formulation of problem (\ref{N-S}), we get
\begin{equation}\label{equality_duality}
\begin{array}{l}
\medskip
\displaystyle
\left\langle F^{\epsilon},\varphi\right\rangle_{W^{-1,C(\delta)}(Q^\epsilon)^3, W^{1,r}_0(Q^\epsilon)^3}\medskip
\\
\displaystyle
=-\epsilon^\gamma\eta\int_{\Omega^\epsilon} Dv^\epsilon: DR^{\epsilon}_r(\varphi)\,dx + \int_{\Omega^\epsilon} f'\cdot (R^{\epsilon}_r(\varphi))'\,dx 
-\int_{\Omega^\epsilon} (v^\epsilon\cdot \nabla)v^\epsilon\,R^{\epsilon}_r(\varphi)\,dx.
\end{array}\end{equation}
By using Lemma \ref{Vel_Estimates_lemma} for fixed $\epsilon$, we will prove that it is a bounded functional on $W^{1,r}_0(Q^\epsilon)$ (see following steps depending on the {\bf TPM}), and so $F^\epsilon\in W^{-1,C(\delta)}(Q^\epsilon)^3$. Moreover, if ${\rm div}\, \varphi=0$, then it holds  $\left\langle F_{\epsilon},\varphi\right\rangle=0\,,$ and so, the DeRham theorem gives the existence of $P^\epsilon$ in $L^{C(\delta)}_0(Q^\epsilon)$ with $F^\epsilon=\nabla P^\epsilon$.\\

\noindent Now, for all $\tilde \varphi\in W^{1,r}_0(\Omega)^3$ with $\tilde v(x', z_3)=v(x',\epsilon z_3)$, by using  (\ref{dilatacion}), we deduce 
$$\begin{array}{l}\displaystyle\langle \nabla_{\epsilon}\tilde P^\epsilon, \tilde \varphi\rangle_{W^{-1,C(\delta)}(\Omega)^3, W^{1,r}_0(\Omega)^3} 
=-\int_{\Omega}\tilde P^\epsilon\,{\rm div}_{\epsilon}\,\tilde \varphi\,dx'dz_3
\medskip\\

\displaystyle =-\epsilon^{-1}\int_{Q^\epsilon}P^\epsilon\,{\rm div}\,\varphi\,dx=\epsilon^{-1}\langle \nabla P^\epsilon, \varphi\rangle_{W^{-1,C(\delta)}(Q^\epsilon)^3, W^{1,r}_0(Q^\epsilon)^3}\,.
\end{array}$$
Taking into account (\ref{equality_duality}), we get
$$\begin{array}{l}\displaystyle\langle \nabla_{\epsilon}\tilde P^\epsilon, \tilde \varphi\rangle_{W^{-1,C(\delta)}(\Omega)^3, W^{1,r}_0(\Omega)^3} \medskip\\
\displaystyle = \epsilon^{-1}\left(
-\epsilon^\gamma\eta\int_{\Omega^\epsilon} Dv^\epsilon: DR^{\epsilon}_r(\varphi)\,dx + \int_{\Omega^\epsilon} f'\cdot (R^{\epsilon}_r(\varphi))'\,dx 
  -\int_{\Omega^\epsilon} (v^\epsilon\cdot \nabla)v^\epsilon\,R^{\epsilon}_r(\varphi)\,dx\right),
\end{array}$$
and after applying (\ref{dilatacion}), we deduce
\begin{equation}\label{extension_1}
\begin{array}{l}\medskip
\displaystyle\langle \nabla_{\epsilon}\tilde P^\epsilon, \tilde \varphi\rangle_{W^{-1,C(\delta)}(\Omega)^3, W^{1,r}_0(\Omega)^3}
\medskip
\\
=\displaystyle- \epsilon^\gamma\eta\int_{\widetilde \Omega^\epsilon} D_\epsilon \tilde v^\epsilon : D_{\epsilon}\tilde R^{\epsilon}_r(\tilde \varphi)\,dx'dz_3+ \int_{\widetilde \Omega^\epsilon} f'\cdot (\tilde R^{\epsilon}_q(\tilde \varphi))'\,dx'dz_3\medskip
\\
\displaystyle 
\quad -\int_{\widetilde \Omega^\epsilon} (\tilde v^\epsilon\cdot \nabla_\epsilon)\tilde v^\epsilon\,\tilde R^{\epsilon}_r(\tilde \varphi)\,dx'dz_3
\end{array}
\end{equation}
with $\tilde R^\epsilon_r(\tilde \varphi)=R^\epsilon_r(\varphi)$ for all $\tilde \varphi \in W^{1,r}_0(\Omega)^3$.\\

\noindent {\it Step 2}. {\it  Estimates of the extended pressure in the case {\bf HTPM}}.  Assume $\delta>1$ and 
we consider an exponent $r\geq 2$ and we denote its conjugate by $C(\delta)\in (1,2]$. Due to the  $L^2$-estimates for velocity, given in (\ref{estimates_v_H}), the value of $r$ must be as close as possible to $2$. Its value will be deduced  when we estimate the intertial term.\\

\noindent Estimate (\ref{estim_restricted_H}) can be rewritten, after applying the change of variables (\ref{dilatacion}), as 
\begin{equation}\label{ext_1} 
\|\tilde R^\epsilon_r\tilde \varphi\|_{L^r(\widetilde\Omega^\epsilon)^3}+ \epsilon^\delta \|D_{\epsilon}\tilde R^\epsilon_r\tilde \varphi\|_{L^r(\widetilde\Omega^\epsilon)^{3\times 3}}\leq  C\left(\|\tilde \varphi\|_{L^r(\Omega)^3} 
+ \epsilon^\delta\|D_{\epsilon}\tilde \varphi\|_{L^r(\Omega)^{3\times 3}}\right),
\end{equation}
and then, as $\epsilon^\delta\ll 1$ and $\epsilon^{-\delta}\gg 1$ because $\delta>1$, we get 
\begin{equation}\label{ext_2}
\|\tilde R^\epsilon_r \tilde \varphi\|_{L^r(\widetilde\Omega^\epsilon)^3}\leq  C \|\tilde \varphi\|_{W^{1,r}_0(\Omega)^3},\quad \|D_\epsilon \tilde R^\epsilon_r\tilde \varphi\|_{L^r(\widetilde\Omega^\epsilon)^{3\times 3}}\leq {C\over \epsilon^\delta}\|\tilde \varphi\|_{W^{1,r}_0(\Omega)^3}.
\end{equation}
Since  $L^r\hookrightarrow L^2 \hookrightarrow L^{C(\delta)}$,  estimates (\ref{estimates_v_H}) and  (\ref{ext_2}), we deduce the estimates for the two first terms of the right-hand side of (\ref{extension_1})
\begin{equation}\label{estim_casoH_1}
\begin{array}{rl}
\medskip
\displaystyle
\left|\epsilon^\gamma\eta  \int_{\widetilde\Omega^\epsilon}  D_\epsilon \tilde v^\epsilon : D_\epsilon\tilde R^{\epsilon}_r(\tilde \varphi)\,dx'dz_3\right|  \leq &\displaystyle C\epsilon^\gamma \|D_{\epsilon} \tilde v^\epsilon\|_{L^{C(\delta)}(\widetilde\Omega^\epsilon)^{3\times 3}}\|D_{\epsilon}\tilde R^\epsilon_r(\tilde \varphi)\|_{L^r(\widetilde \Omega^\epsilon)^{3\times 3}}\\
\medskip
\leq &\displaystyle C\epsilon^\gamma \|D_{\epsilon}\tilde v^\epsilon\|_{L^{2}(\widetilde\Omega^\epsilon)^{3\times 3}}\|D_{\epsilon}\tilde R^\epsilon_r(\tilde \varphi)\|_{L^r(\widetilde \Omega^\epsilon)^{3\times 3}}
\\
\medskip
 \leq &\displaystyle  C \|\tilde \varphi\|_{W^{1,r}_0(\Omega)^3},
\end{array}
\end{equation}
\begin{equation}\label{estim_casoH_2}
\begin{array}{rl}
\medskip
\displaystyle
 \left|\int_{\widetilde\Omega^\epsilon}f'\cdot (\tilde R^\epsilon_r (\tilde \varphi))' \,dx'dz_3\right|&\leq C\|\tilde R^\epsilon_r (\tilde \varphi) \|_{L^r(\widetilde\Omega^\epsilon)^3}\leq C\|\tilde \varphi\|_{W^{1,r}_0(\Omega)^3}\,.
\end{array}
\end{equation}
To finish, we must estimate the inertial term. Observe that the inertial term can be written as 
\begin{equation}\label{inertial_caseH_3}
\begin{array}{l}
\displaystyle \int_{\widetilde \Omega^\epsilon} (\tilde v^\epsilon\cdot \nabla_{\epsilon}) \tilde v^\epsilon\,\tilde R^\epsilon_r( \tilde \varphi)\, dx'dz_3\medskip
\\
= \displaystyle -\int_{\widetilde \Omega^\epsilon} \tilde v^\epsilon\tilde \otimes \tilde v^\epsilon:D_{x'} \tilde R^\epsilon_r (\tilde \varphi)\,dx'dz_3
\medskip
\\
\quad \displaystyle +{1\over \epsilon}\left(\int_{\widetilde \Omega^\epsilon}\partial_{z_3} \tilde v^{\epsilon}_3\tilde  v^\epsilon \tilde R^\epsilon_r (\tilde \varphi) \,dx'dz_3+ \int_{\widetilde \Omega^\epsilon}\tilde v^{\epsilon}_3\partial_{z_3} \tilde v^\epsilon\, \tilde R^\epsilon_r(\tilde\varphi)\,dx'dz_3\right),
\end{array}
\end{equation}
where $(\varphi \tilde \otimes \psi)_{ij}=\varphi_i \psi_j$, $i=1,2$, $j=1,2,3$. \\

\noindent We will prove that for the values of $C(\delta)$ given in (\ref{C_case_H}) in terms of $\delta$ and $\gamma$, the inertial term will vanish in the limit when  $\epsilon$ tends to zero. To do this, it is enough to prove that the inertial term (\ref{inertial_caseH_3}) satisfies
\begin{equation}\label{inertial_alpha}\left|\int_{\widetilde\Omega^\epsilon}(\tilde v^\epsilon\cdot \nabla_\epsilon)\tilde v^\epsilon \tilde R^\epsilon_r (\tilde \varphi)\,dx'dz_3\right|\leq C\epsilon^{\alpha}\|\tilde \varphi\|_{W^{1,r}_0(\Omega)^3},\quad \forall\, \tilde \varphi \in W^{1,r}_0(\Omega)^3.
\end{equation}
where $\alpha\geq 0$. Below, we estimate each term of the right-hand side of (\ref{inertial_caseH_3}).

\begin{itemize}
\item {\it First term in the second member of (\ref{inertial_caseH_3})}. From H\"older's inequality and  the second estimate of  (\ref{ext_2}), we have
\begin{equation}\label{interpola_1_caso2}
\begin{array}{rl}
\displaystyle \left|\int_{\widetilde \Omega^\epsilon} \tilde v^\epsilon\tilde \otimes \tilde v^\epsilon:D_{x'} \tilde R^\epsilon_r( \tilde \varphi)\,dx'dz_3\right| 
&\displaystyle \leq \|\tilde v^\epsilon\|_{L^{r'}(\widetilde\Omega^\epsilon)^3}^2\|D_{x'}\tilde R^\epsilon_r(\tilde \varphi)\|_{L^r(\widetilde\Omega^\epsilon)^{3\times 2}}\medskip\\
&\displaystyle \leq C\epsilon^{-\delta}  \|\tilde v^\epsilon\|_{L^{r'}(\widetilde\Omega^\epsilon)^3}^2\|\tilde \varphi\|_{W^{1,r}_0(\Omega)^3},
\end{array}
\end{equation}
where $r$ and $r'$ satisfy
\begin{equation}\label{ineq1}
{2\over r'}+{1\over r}\leq 1.
\end{equation}
We need to estimate $ \|\tilde v^\epsilon\|_{L^{r'}(\widetilde\Omega^\epsilon)^3}$ taking into account estimates (\ref{estimates_v_H}). For that reason, our objetive is to consider $r$ as near as possible to $2$. Thus,  considering an interpolation between $L^2(\widetilde \Omega^\epsilon)$ and $H^1_0(\widetilde \Omega^\epsilon)$ introducing the parameter $\theta\in [0,1]$ such that 
\begin{equation}\label{inter2_6}{1\over r'}={1\over 2}\theta+{1\over 6}(1-\theta),
\end{equation}
we deduce
\begin{equation*} \|\tilde v^\epsilon\|_{L^{r'}(\widetilde\Omega^\epsilon)^3}\leq \|\tilde v^\epsilon\|^{\theta}_{L^2(\widetilde\Omega^\epsilon)^3}\|\tilde v^\epsilon\|^{1-\theta}_{L^{6}(\widetilde\Omega^\epsilon)^3}.
\end{equation*}
From  $H^1_0(\widetilde \Omega^\epsilon) \hookrightarrow L^{6}(\widetilde\Omega^\epsilon)$ and  estimates (\ref{estimates_v_H}), we get
\begin{equation}\label{interpola_2_caso2} 
\begin{array}{rl}
\|\tilde v^\epsilon\|_{L^{r'}(\widetilde\Omega^\epsilon)^3}& \displaystyle \leq \|\tilde v^\epsilon\|^{\theta}_{L^2(\widetilde\Omega^\epsilon)^3}\| D \tilde v^\epsilon\|^{1-\theta}_{L^2(\widetilde\Omega^\epsilon)^{3\times 3}}\medskip\\
&\displaystyle \leq \|\tilde v^\epsilon\|^{\theta}_{L^2(\widetilde\Omega^\epsilon)^3}\| D_\epsilon \tilde v^\epsilon\|^{1-\theta}_{L^2(\widetilde\Omega^\epsilon)^{3\times 3}}
\leq C\epsilon^{\theta(2\delta-\gamma)+(1-\theta)(\delta-\gamma)}.
\end{array}
\end{equation}
Taking into account (\ref{interpola_2_caso2}) in (\ref{interpola_1_caso2}), we deduce
\begin{equation}\label{interpola_3_caso2}
\begin{array}{l}
\displaystyle \left|\int_{\widetilde \Omega^\epsilon} \tilde v^\epsilon\tilde \otimes \tilde v^\epsilon:D_{x'} \tilde R^\epsilon_r( \tilde \varphi)\,dx'dz_3\right| 
\leq  C\epsilon^{2\delta\left(\theta-{\gamma\over \delta}+{1\over 2}\right)} \|\tilde \varphi\|_{W_0^{1,r}( \Omega )^{3}}.
\end{array}
\end{equation}
Thus, we must choose $\theta$ satisfying $\theta-\gamma/\delta+1/2\geq 0$, i.e.
\begin{equation}\label{intertial_case1_5}
\theta\geq\theta_0= \max\left\{0,{\gamma\over \delta}-{1\over2}\right\}.
\end{equation}
From the inequality (\ref{ineq1}) and the equation (\ref{inter2_6}), we deduce  that 
\begin{equation}\label{q_bounded}
r\ge {3\over 2(1-\theta)},
\end{equation}
and from (\ref{intertial_case1_5}) and (\ref{q_bounded}), we get
\begin{equation}\label{intertial_case1_6}
r\ge {3\delta\over 3\delta-2\gamma}.
\end{equation}
Observe that if we consider $r=2$ in (\ref{intertial_case1_6}), we get that $\gamma\leq (3\delta)/ 4$. For $\gamma>(3\delta)/4$,  we consider (\ref{intertial_case1_6}) taking into account that $\gamma/\delta-1/2\leq \theta \leq 1$, i.e.  we consider (\ref{intertial_case1_6}) for $(3\delta)/4<\gamma\leq (3\delta)/2$. Therefore, we obtain
\begin{equation}\label{case_1_q_term1} \left\{\begin{array}{lcl}
r=2& \hbox{for}& \gamma \leq \displaystyle{3\delta\over 4},\\
\\
\displaystyle r\geq{3\delta\over 3\delta-2\gamma}& \hbox{for}&\displaystyle {3\delta\over 4}<\gamma \leq {3\delta\over 2}.
\end{array}\right.
\end{equation}

\item {\it Second and third terms in the second member of (\ref{inertial_caseH_3})}. By using H${\rm \ddot{o}}$lder's inequality, the second and third terms are estimated by 
\begin{equation}\label{estimate_(b)}
\begin{array}{l}
\displaystyle \left|\epsilon^{-1}\left(\int_{\widetilde \Omega^\epsilon}\partial_{z_3} \tilde v^{\epsilon}_3\tilde  v^\epsilon \tilde R^\epsilon_r (\tilde \varphi) \,dx'dz_3+ \int_{\widetilde \Omega^\epsilon}\tilde v^{\epsilon}_3\partial_{z_3} \tilde v^\epsilon\, \tilde R^\epsilon_r (\tilde \varphi)\,dx'dz_3\right)\right|
\medskip
\\
\displaystyle \leq C \epsilon^{-1}\|\partial_{z_3}\tilde v^\epsilon\|_{L^2(\widetilde\Omega^\epsilon)^3}\|\tilde v^\epsilon\|_{L^{r'}(\widetilde\Omega^\epsilon)^3}\| R^\epsilon_r \tilde v\|_{L^{r_1}(\widetilde\Omega^\epsilon)^3}, 
\end{array}
\end{equation}
where $r'$ and $r_1$ satisfy
\begin{equation}\label{ineq2}
{1\over 2}+{1\over r'}+{1\over r_1}\leq 1.
\end{equation}
Interpolating between $L^2(\widetilde \Omega^\epsilon)$ and $H^1_0(\widetilde \Omega^\epsilon)$ with  $\theta\in [0,1]$, we have again relation (\ref{inter2_6}) and estimate (\ref{interpola_2_caso2}) of $\tilde v^\epsilon$ in $L^{r'}$, i.e. 
$${1\over r'}={1\over 2}\theta+{1\over 6}(1-\theta)\quad \hbox{and}\quad  \|\tilde v^\epsilon\|_{L^{r'}(\widetilde\Omega^\epsilon)^3}
\leq C\epsilon^{\theta(2\delta-\gamma)+(1-\theta)(\delta-\gamma)}.
$$

Taking into account this and  the   second estimate in (\ref{estimates_v_H}), then   (\ref{estimate_(b)}) satisfy
\begin{equation}\label{estimate_(b)2}
\begin{array}{l}
\displaystyle \left|\epsilon^{-1}\left(\int_{\widetilde \Omega^\epsilon}\partial_{z_3} \tilde v^{\epsilon}_3\tilde  v^\epsilon \tilde R^\epsilon_r (\tilde \varphi) \,dx'dz_3+ \int_{\widetilde \Omega^\epsilon}\tilde v^{\epsilon}_3\partial_{z_3} \tilde v^\epsilon\, \tilde R^\epsilon_r (\tilde \varphi)\,dx'dz_3\right)\right|\medskip
\\
\displaystyle \leq C\epsilon^{\theta\delta+2\delta-2\gamma}\| R^\epsilon_r(\tilde \varphi)\|_{L^{r_1}(\widetilde\Omega^\epsilon)^3}.
\end{array}
\end{equation}
Finally, it remains estimate $\| R^\epsilon_r(\tilde \varphi)\|_{L^{r_1}}$. To do this, we use the Sobolev-Gagliardo-Nirenberg theorem, which states that if $1\leq r<3$, then we have $W^{1,r} \hookrightarrow L^{r^*}$ with $1/r^*=1/r-1/3$. As our goal is to have $\tilde \varphi\in W^{1,r}(\Omega)^3$ with $r\ge 2$ as close as possible to $2$, we choose $r_1={3r\over 3-r}$ which satisfies (\ref{ineq2}).\\

 Then, by using the Sobolev-Gagliardo-Nirenberg theorem and the second estimate given in (\ref{ext_2}) imply that
$$\begin{array}{rl}
\displaystyle \| R^\epsilon_r( \tilde \varphi)\|_{L^{r_1}(\widetilde\Omega^\epsilon)^3} &\displaystyle 
\leq C \|D R^\epsilon_r (\tilde \varphi)\|_{L^{r}(\widetilde\Omega^\epsilon)^{3\times 3}}\medskip\\
&\displaystyle \leq C \|D_\epsilon R^\epsilon_r (\tilde \varphi)\|_{L^{r}(\widetilde\Omega^\epsilon)^{3\times 3}}\leq C\epsilon^{-\delta}\|\tilde \varphi\|_{W^{1,r}_0(\Omega)^3},
\end{array}$$
and so, (\ref{estimate_(b)2}) reads as follow
\begin{equation}\label{estimate_(b)3}
\begin{array}{l}\displaystyle
\left|\epsilon^{-1}\left(\int_{\widetilde \Omega^\epsilon}\partial_{z_3} \tilde v^{\epsilon}_3\tilde  v^\epsilon \tilde R^\epsilon_r (\tilde \varphi) \,dx'dz_3+ \int_{\widetilde \Omega^\epsilon}\tilde v^{\epsilon}_3\partial_{z_3} \tilde v^\epsilon\, \tilde R^\epsilon_r (\tilde \varphi)\,dx'dz_3\right)\right|\medskip
\\
\displaystyle \leq  C\epsilon^{\delta(\theta+1-{2\over \delta}\gamma)}\|\tilde \varphi\|_{W^{1,r}_0(\Omega)^3}.
\end{array}
\end{equation}
Thus,  $\theta$ has to satisfy $\theta-2\gamma/\delta+1\geq 0$, i.e.
\begin{equation}\label{intertial_case1_52}
\theta\geq \max\left\{0, {2\over \delta}\gamma-1\right\}.
\end{equation}
The inequality (\ref{ineq2}) with $r_1={3r\over 3-r}$ together with the equation (\ref{inter2_6}) imply that 
\begin{equation}\label{relation_min_q}
\theta\leq \min\left\{1,2-{3\over r}\right\}.
\end{equation}
 Then,  since $ 2\gamma/ \delta -1\leq \theta\leq 2-{3/r}$, then  we deduce (\ref{intertial_case1_6}).\\

Assuming $r=2$ in (\ref{intertial_case1_6}), we deduce that  $\gamma\leq 3\delta/4$. For $\gamma>3\delta/4$, we consider (\ref{intertial_case1_6}) taking into account that $2\gamma/ \delta -1\leq \theta \leq 1$, then $\gamma\leq \delta$ and so, we consider (\ref{intertial_case1_6}) for $3\delta/4<\gamma\leq \delta$. Therefore, we finally get
\begin{equation}\label{case_1_q_term2} \left\{\begin{array}{lcl}
r=2& \hbox{for}& \gamma \leq \displaystyle  {3\delta\over 4},\\
\\
\displaystyle r\geq {3\delta\over 3\delta-2\gamma}& \hbox{for}&\displaystyle{3\delta\over 4}<\gamma\leq \delta.
\end{array}\right.
\end{equation}
\end{itemize}
Then, considering the weaker condition between (\ref{case_1_q_term1}) and (\ref{case_1_q_term2}), we have that the desired estimate (\ref{inertial_alpha}) holds by considering $\delta>1$ and $r$ as follows
\begin{equation}\label{case_1_q_termfinal} \left\{\begin{array}{lcl}
r=2& \hbox{for}& \gamma \leq \displaystyle {3\delta\over 4},\\
\\
\displaystyle  r\geq{3\delta\over 3\delta-2\gamma}& \hbox{for}&\displaystyle {3\delta\over 4}<\gamma \leq  \delta.
\end{array}\right.
\end{equation}
To finish, denoting $C(\delta)$ the conjugate of $r$, defined by (\ref{C_case_H}), and taking into account estimates (\ref{estim_casoH_1}), (\ref{estim_casoH_2}) and (\ref{inertial_alpha})  in (\ref{extension_1}), we get the  estimate for $\nabla_\epsilon\tilde P^\epsilon$ given in (\ref{esti_p_caseH}).  Taking into account the Ne${\check{\rm c}}$as inequality, there exists a representative $\tilde P^\epsilon\in L^{C(\delta)}_0(\Omega)$  such that
$$\|\tilde P^\epsilon\|_{L^{C(\delta)}(\Omega)}\leq C\|\nabla\tilde P^\epsilon\|_{W^{-1,C(\delta)}(\Omega)^3}\leq C\|\nabla_{\epsilon}\tilde P^\epsilon\|_{W^{-1,C(\delta)}(\Omega)^3},$$
which gives the first estimate in (\ref{esti_p_caseH}).
\\

\noindent {\it Step 3}. {\it  Estimates of the extended pressure in the case {\bf PTPM}}.  Assume $\delta=1$. As the previous case, we have to consider an exponent $r\geq 2$. Due to the  $L^2$-estimates for velocity, given in (\ref{estimates_v_P}), the value of $r$ must be as close as possible to $2$. Its value can be deduced after estimating the intertial term.  In this case, the dilated restricted operator satisfies the following estimates
\begin{equation}\label{ext_1_P}\begin{array}{l}
\medskip\displaystyle
\|\tilde R^\epsilon_r\tilde \varphi\|_{L^r(\widetilde\Omega^\epsilon)^3}+ \epsilon \|D_{\epsilon}\tilde R^\epsilon_r\tilde \varphi\|_{L^r(\widetilde\Omega^\epsilon)^{3\times 3}}\leq  C\left(\|\tilde \varphi\|_{L^r(\Omega)^3} 
+ \epsilon \|D_{\epsilon}\tilde \varphi\|_{L^r(\Omega)^{3\times 3}}\right),
\end{array} 
\end{equation}
and then, as $\epsilon\ll 1$, it holds
\begin{equation}\label{ext_2_P}
\|\tilde R^\epsilon_r \tilde \varphi\|_{L^r(\widetilde\Omega^\epsilon)^3}\leq  C \|\tilde \varphi\|_{W^{1,r}_0(\Omega)^3},\quad \|D_\epsilon \tilde R^\epsilon_r\tilde \varphi\|_{L^r(\widetilde\Omega^\epsilon)^{3\times 3}}\leq {C\over \epsilon}\|\tilde \varphi\|_{W^{1,r}_0(\Omega)^3}.
\end{equation}
The proof follows the lines of {\it Step 2} by taking into account estimates of the velocity (\ref{estimates_v_P}), estimates of the restricted operator (\ref{ext_2_P}) and just changing the value of $\delta$ by $1$, so we omit it. As a result, we can deduce that for $r$ satisfying (\ref{case_1_q_termfinal}) imposing $\delta=1$, it holds
\begin{equation}\label{estim_casoP_1}
\begin{array}{rl}
\medskip
\displaystyle
\left|\epsilon^\gamma\eta  \int_{\widetilde\Omega^\epsilon}  D_\epsilon \tilde v^\epsilon : D_\epsilon\tilde R^{\epsilon}_r(\tilde \varphi)\,dx'dz_3\right|  \leq &  C \|\tilde \varphi\|_{W^{1,r}_0(\Omega)^3},
\end{array}
\end{equation}
\begin{equation}\label{estim_casoP_2}
\begin{array}{rl}
\medskip
\displaystyle
 \left|\int_{\widetilde\Omega^\epsilon}f'\cdot (\tilde R^\epsilon_r (\tilde \varphi))' \,dx'dz_3\right|&\leq  C\|\tilde \varphi\|_{W^{1,r}_0(\Omega)^3}\,,
\end{array}
\end{equation}
and 
\begin{equation}\label{inertial_alphaP}\left|\int_{\widetilde\Omega^\epsilon}(\tilde v^\epsilon\cdot \nabla_\epsilon)\tilde v^\epsilon \tilde R^\epsilon_r (\tilde \varphi)\,dx'dz_3\right|\leq C\epsilon^{\alpha}\|\tilde \varphi\|_{W^{1,r}_0(\Omega)^3},\quad \alpha\geq 0.
\end{equation}
Then, considering $C(\delta)$ given by (\ref{C_case_H}), it holds that $C(1)$ is the conjugate of $r$ in this case. Taking into account estimates (\ref{estim_casoP_1}), (\ref{estim_casoP_2}) and (\ref{inertial_alphaP})  in (\ref{extension_1}), we get the  estimate for $\nabla_\epsilon\tilde P^\epsilon$ given in (\ref{esti_p_caseP}).  Taking into account the Ne${\check{\rm c}}$as inequality, there exists a representative $\tilde P^\epsilon\in L^{C(1)}_0(\Omega)$  such that
$$\|\tilde P^\epsilon\|_{L^{C(1)}(\Omega)}\leq C\|\nabla\tilde P^\epsilon\|_{W^{-1,C(1)}(\Omega)^3}\leq C\|\nabla_{\epsilon}\tilde P^\epsilon\|_{W^{-1,C(1)}(\Omega)^3},$$
which gives the first estimate in (\ref{esti_p_caseP}).\\

\noindent {\it Step 4}. {\it  Estimates of the extended pressure in the case {\bf VTPM}}.  Assume $\delta\in (0,1)$. The velocity estimates in this case are given by (\ref{estimates_v_V}) and   the dilated restricted operator satisfies the following estimates
\begin{equation}\label{ext_1_V}\begin{array}{l}
\medskip\displaystyle
\|\tilde R^\epsilon_r\tilde \varphi\|_{L^r(\widetilde\Omega^\epsilon)^3}+ \epsilon \|D_{\epsilon}\tilde R^\epsilon_r\tilde \varphi\|_{L^r(\widetilde\Omega^\epsilon)^{3\times 3}}\leq  C\left(\|\tilde \varphi\|_{L^r(\Omega)^3} 
+ \epsilon^\delta \|D_{\epsilon}\tilde \varphi\|_{L^r(\Omega)^{3\times 3}}\right),
\end{array} 
\end{equation}
and then, as $\epsilon^\delta\ll 1$ and $\epsilon^{\delta-1}\ll\epsilon^{-1}$ because $\delta\in (0,1)$, it holds
\begin{equation}\label{ext_2_V}
\|\tilde R^\epsilon_r \tilde \varphi\|_{L^r(\widetilde\Omega^\epsilon)^3}\leq  C \|\tilde \varphi\|_{W^{1,r}_0(\Omega)^3},\quad \|D_\epsilon \tilde R^\epsilon_r\tilde \varphi\|_{L^r(\widetilde\Omega^\epsilon)^{3\times 3}}\leq {C\over \epsilon}\|\tilde \varphi\|_{W^{1,r}_0(\Omega)^3}.
\end{equation}
Since the velocity estimates and the restriced estimates are the same as in the case {\bf PTPM}, respectively, the proof follows the lines of {\it Step 2} with $\delta=1$, so we omit it. Thus,  we get the estimate for $\nabla_\epsilon\tilde P^\epsilon$ given in (\ref{esti_p_caseV}) and so, taking into account the Ne${\check{\rm c}}$as inequality, there exists a representative $\tilde P^\epsilon\in L^{C(1)}_0(\Omega)$  such that
$$\|\tilde P^\epsilon\|_{L^{C(1)}(\Omega)}\leq C\|\nabla\tilde P^\epsilon\|_{W^{-1,C(1)}(\Omega)^3}\leq C\|\nabla_{\epsilon}\tilde P^\epsilon\|_{W^{-1,C(1)}(\Omega)^3},$$
which gives the first estimate in (\ref{esti_p_caseV}).
\par\hfill$\square$

\subsection{Unfolded functions}\label{Subsec:unfolded}
Following \cite{Anguiano_SuarezGrau, Anguiano_SuarezGrau2},  to have information about the microstructure of the domain $\widetilde\Omega^\epsilon$, we introduce the adaptation of the unfolding method from \cite{CDG}. This adaptation divides the domain $\widetilde\Omega^\epsilon$ in rectangular parallelepipeds of lateral length $\epsilon^\delta$ and vertical length $1$. Thus, given $(\tilde \varphi^{\epsilon},  \tilde \psi^\epsilon) \in H^1_0(\Omega)^3\times L^{r}_0(\Omega)$, $1<r<+\infty$, we define $(\hat{\varphi}_{\epsilon},  \hat \psi^\epsilon)$ by
\begin{equation}\label{phihat}
\hat{\varphi}^{\epsilon}(x^{\prime},z)=\tilde{\varphi}^{\epsilon}\left( {\epsilon}^\delta\kappa\left(\frac{x^{\prime}}{{\epsilon^\delta }} \right)+{\epsilon}^\delta z^{\prime},z_3 \right),\quad 
\hat{\psi}^{\epsilon}(x^{\prime},z)=\tilde{\psi}^{\epsilon}\left( {\epsilon}^\delta \kappa\left(\frac{x^{\prime}}{{\epsilon^\delta }} \right)+{\epsilon}^\delta z^{\prime},z_3 \right),
\end{equation}
  a.e. $(x^{\prime},z)\in \omega\times Z$, assuming $\tilde \varphi^\epsilon$ and $\tilde \psi^\epsilon$ are extended by zero outside $\omega$, where the function $\kappa:\mathbb{R}^2\to \mathbb{Z}^2$ is defined by 
$$\kappa(x')=k'\Longleftrightarrow x'\in Z'_{k',1},\quad\forall\,k'\in\mathbb{Z}^2.$$
\begin{remark}\label{remarkCV_1}
The function $\kappa$ is well defined up to a set of zero measure in $\mathbb{R}^2$ (the set $\cup_{k'\in \mathbb{Z}^2}\partial Z'_{k',1}$). Moreover, for every $\epsilon>0$, we have
$$\kappa\left({x'\over \epsilon^\delta }\right)=k'\Longleftrightarrow x'\in Z'_{k',\epsilon^\delta }.$$
\end{remark}
\noindent We have the following estimates relating  $(\hat \varphi^\epsilon,\hat \psi^\epsilon)$ and $(\tilde \varphi^\epsilon, \tilde \psi^\epsilon)$, see the proof of \cite[Lemma 4.9]{Anguiano_SuarezGrau} for more details. \\

\begin{lemma}\label{estimates_relation} The sequence $(\hat \varphi^\epsilon, \hat \psi^\epsilon)$ defined by (\ref{phihat}) satisfies
\begin{equation}\label{estimates_relation_1}
\begin{array}{c}
\displaystyle
\medskip\|\hat \varphi^\epsilon\|_{L^2(\omega\times Z)^3}\leq \|\tilde \varphi^\epsilon\|_{L^2(\Omega)^3},
\medskip\\
\medskip
\displaystyle \|D_{z'} \hat \varphi^\epsilon\|_{L^2(\omega\times Z)^{3\times 2}}\leq \epsilon^\delta \|D_{x'}\tilde \varphi^\epsilon\|_{L^2(\Omega)^{3\times 2}},\quad \|\partial_{z_3} \hat \varphi^\epsilon\|_{L^2(\omega\times Z)^{3 }}\leq   \|\partial_{z_3}\tilde \varphi^\epsilon\|_{L^2(\Omega)^{3}},
\medskip
\\
\|\hat \psi^\epsilon\|_{L^{r}(\omega\times Z)}\leq \|\tilde \psi^\epsilon\|_{L^{r}(\Omega)}.
\end{array}
\end{equation}
\end{lemma}
\begin{definition}[Unfolded velocity and pressure]\label{Definition_unfolded} Consider $\delta>0$ and the function $C(\delta)$ given by (\ref{C_case_H}). We define the unfolded velocity and pressure $(\hat v^\epsilon, \hat P^\epsilon)$ from  $(\tilde v^\epsilon, \tilde P^\epsilon)$ depending on   the type of {\bf TPM}:
\begin{itemize}
\item[--] In the case {\bf HTPM} ($\delta>1$) and $\gamma\leq \delta$,  from $(\tilde v^\epsilon, \tilde P^\epsilon)\in H^1_0(\Omega)^3\times L^{C(\delta)}_0(\Omega)$,  we define $(\hat v^\epsilon, \hat P^\epsilon)$ by  using (\ref{phihat}) with $\tilde \varphi^\epsilon=\tilde v^\epsilon$ and $\tilde \psi^\epsilon=\tilde P^\epsilon$.\\

\item[--] In the case {\bf PTPM} ($\delta=1$) and $\gamma\leq 1$,  from $(\tilde v^\epsilon, \tilde P^\epsilon)\in H^1_0(\Omega)^3\times L^{C(1)}_0(\Omega)$, we define $(\hat v^\epsilon, \hat P^\epsilon)$ by  using (\ref{phihat}) with $\tilde \varphi^\epsilon=\tilde v^\epsilon$ and $\tilde \psi^\epsilon=\tilde P^\epsilon$.\\

\item[--] In the case {\bf VTPM} ($0<\delta<1$) and   $\gamma\leq 1$,  from $(\tilde v^\epsilon, \tilde P^\epsilon)\in H^1_0(\Omega)^3\times L^{C(1)}_0(\Omega)$, we define $(\hat v^\epsilon, \hat P^\epsilon)$ by  using (\ref{phihat}) with $\tilde \varphi^\epsilon=\tilde v^\epsilon$ and $\tilde \psi^\epsilon=\tilde P^\epsilon$.\\
\end{itemize}
\end{definition}
\begin{remark}\label{remarkCV2}
For $k^{\prime}\in \mathcal{K}_{\epsilon}$, the restrictions of $(\hat{v}^{\epsilon},   \hat P^\epsilon)$  to $Z^{\prime}_{k^{\prime},{\epsilon^\delta}}\times Z$ does not depend on $x^{\prime}$, whereas as a function of $z$ it is obtained from $(\tilde{v}^{\epsilon},  \tilde{P}^{\epsilon})$ by using the change of variables $\displaystyle z^{\prime}=\frac{x^{\prime}- {\epsilon}^\delta k^{\prime}}{{\epsilon^\delta}}$,
which transforms $Z_{k^{\prime}, {\epsilon^\delta}}$ into $Z$.\\
\end{remark}

\noindent From the estimates of the extension of the velocity given in Lemma \ref{Vel_Estimates_lemma}, the estimates of the extension of the pressure given in Theorem \ref{mainthmEstimates} and the estimates given in Lemma \ref{estimates_relation}, we deduce the following estimates for  $(\hat v^\epsilon,\hat P^\epsilon)$ in each {\bf TMP}.\\

\begin{lemma}[Estimates of the unfolded functions]\label{estimates_hat} Consider $\gamma\in \mathbb{R}$, $\delta>0$ and $C(\delta)$ given by (\ref{C_case_H}). Depending on the type of {\bf TPM}, we have the following estimates:
 \begin{itemize}
 \item[(i)] In the case {\bf HTPM} ($\delta>1$) and $\gamma\leq \delta$,  there exists a constant $C>0$, independent of $\epsilon$, such that
\begin{equation}  \label{estim_u_hat}
\begin{array}{c}\displaystyle
 \|\hat v^\epsilon\|_{L^2(\omega\times Z)^3}\leq C\epsilon^{2\delta-\gamma},\medskip\\
 \displaystyle  
 \|D_{z'}\hat v^\epsilon\|_{L^2(\omega\times Z)^{3\times 2}}\leq C\epsilon^{2\delta-\gamma},\quad  \|\partial_{z_3}\hat v^\epsilon\|_{L^2(\omega\times Z)^{3}}\leq C\epsilon^{\delta+1-\gamma},
 \end{array}
\end{equation}
\begin{equation}\label{estim_P_hat}
 \|\hat P^\epsilon\|_{L^{C(\delta)}(\omega\times Z)}\leq C.
\end{equation}

\item[(ii)] In the case {\bf PTPM} ($\delta=1$) and $\gamma\leq 1$,  there exists a constant $C>0$, independent of $\epsilon$, such that
\begin{equation}  \label{estim_u_hat}
 \|\hat v^\epsilon\|_{L^2(\omega\times Z)^3}\leq C\epsilon^{2-\gamma}, \quad
 \|D_{z}\hat v^\epsilon\|_{L^2(\omega\times Z)^{3\times 3}}\leq C\epsilon^{2-\gamma},
\end{equation}
\begin{equation}\label{estim_P_hat}
 \|\hat P^\epsilon\|_{L^{C(1)}(\omega\times Z)}\leq C.
\end{equation}

\item[(iii)] In the case {\bf VTPM} ($0<\delta<1$) and $\gamma\leq 1$,  there exists a constant $C>0$, independent of $\epsilon$, such that
\begin{equation}  \label{estim_u_hat}\begin{array}{c}\displaystyle
 \|\hat v^\epsilon\|_{L^2(\omega\times Z)^3}\leq C\epsilon^{2-\gamma}, \medskip\\
 \displaystyle
 \|D_{z'}\hat v^\epsilon\|_{L^2(\omega\times Z)^{3\times 2}}\leq C\epsilon^{\delta+1-\gamma},\quad  \|\partial_{z_3}\hat v^\epsilon\|_{L^2(\omega\times Z)^{3}}\leq C\epsilon^{2-\gamma},
 \end{array}
\end{equation}
\begin{equation}\label{estim_P_hat}
 \|\hat P^\epsilon\|_{L^{C(1)}(\omega\times Z)}\leq C.
\end{equation}
 \end{itemize}
 \end{lemma}

\subsection{Variational formulation of the unfolded functions}\label{Subsec:varfun}
Below, we derive the variational formulation of the unfolded functions, which  will be useful to pass to the limit later. We will give more in detail the inertial term, which is the novelty in this paper. For more information concerning the derivation, see \cite[Theorem 6.1]{Anguiano_SuarezGrau}.\\
 
We consider  $\varphi(x',z)\in \mathcal{D}(\omega; C^\infty_{\#}(Z)^3)$ with $\varphi(x',z)=0$ in $\omega\times S$. Considering  $\varphi^\epsilon(x',z_3)=\varphi(x',x'/\epsilon^\delta,z_3)$ as test function in (\ref{N-S-d}), integrating by parts and taking into account the extension of $\tilde v^\epsilon$ and $\tilde P^\epsilon$, we have
\begin{equation}\label{variational_formulation_tilde}
\begin{array}{l}
\displaystyle
\epsilon^\gamma\eta \int_{ \Omega}D_{x'} \tilde v^\epsilon :\left( D_{x'}\varphi + \epsilon^{-\delta}D_{z'}\varphi \right)\,dx'dz_3+\epsilon^{\gamma-2}\eta \int_{ \Omega} \partial_{z_3}\tilde v^\epsilon\cdot \partial_{z_3}\varphi\,dx'dz_3\medskip\\
\displaystyle+\int_{\Omega}\left(\tilde v^\epsilon\cdot \nabla_{\epsilon}\right)\tilde v^\epsilon\,\varphi\,dx'dz_3
-\int_{\Omega}\tilde P^\epsilon \left({\rm div}_{x'}\varphi'+\epsilon^{-\delta}{\rm div}_{z'}\varphi'+\epsilon^{-1}\partial_{z_3}\varphi_3\right)\,dx'dz_3\medskip\\
\displaystyle=\int_{\Omega} f'\cdot \varphi'\,dx'dz_3.
\end{array}
\end{equation}
 The inertial term can be rewritten as 
\begin{equation}\label{inertial_variational_formulation}
\begin{array}{rl}
\displaystyle \int_{ \Omega} \left(\tilde v^\epsilon\cdot \nabla_{\epsilon}\right)\tilde v^\epsilon\,\varphi\,dx'dz_3=&
\displaystyle 
-\int_{ \Omega} \tilde v^\epsilon\tilde \otimes \tilde v^\epsilon:(D_{x'}  \varphi+\epsilon^{-\delta}D_{z'}\varphi)\,dx'dz_3
\medskip
\\ 
&\displaystyle +\epsilon^{-1}\left(\int_{ \Omega}\partial_{z_3} \tilde v^{\epsilon}_3\tilde  v^\epsilon   \varphi \,dx'dz_3
+ \int_{ \Omega}\tilde v^{\epsilon}_3\partial_{z_3} \tilde v^\epsilon  \varphi\,dx'dz_3\right)\,.
\end{array}
\end{equation}
Applying the unfolding change of variables (see Remark \ref{remarkCV2})  in (\ref{inertial_variational_formulation}), we get
\begin{equation}\label{inertial_variational_formulation_unfolded}
\begin{array}{l}
\displaystyle \int_{ \Omega} \left(\tilde v^\epsilon\cdot \nabla_{\epsilon}\right)\tilde v^\epsilon\varphi\,dx'dz_3
\medskip \\  
\displaystyle=   -\epsilon^{-\delta}\int_{ \omega\times Z} \hat v^\epsilon\tilde \otimes \hat v^\epsilon:D_{z'}  \varphi\,dx'dz\medskip\\
\displaystyle\quad  +\epsilon^{-1}\left(\int_{ \omega\times Z}\partial_{z_3} \hat v^{\epsilon}_3\hat  v^\epsilon  \varphi \,dx'dz+ \int_{ \omega\times Z}\hat v^{\epsilon}_3\partial_{z_3} \hat v^\epsilon\,  \varphi\,dx'dz\right) +O_\epsilon,
\end{array}
\end{equation}
where $O_\epsilon$ tends to zero and can change from line to line.\\

\noindent With this and  by applying the unfolding change of variables  to the rest of the terms in (\ref{variational_formulation_tilde}), we get
\begin{equation}\label{form_var_hat}\begin{array}{l}
\displaystyle
\epsilon^{\gamma-2\delta}\eta\int_{ \omega\times Z}D_{z'}\hat v^\epsilon:D_{z'}\varphi\,dx'dz+
\epsilon^{\gamma-2}\eta\int_{ \omega\times Z}\partial_{z_3}\hat v^\epsilon:\partial_{z_3}\varphi\,dx'dz
\medskip\\
\displaystyle
-\epsilon^{-\delta}\int_{ \omega\times Z} \hat v^\epsilon\tilde \otimes \hat v^\epsilon:D_{z'}  \varphi\,dx'dz
\displaystyle \medskip\\
\displaystyle+\epsilon^{-1}\left(\int_{ \omega\times Z}\partial_{z_3} \hat v^{\epsilon}_3\hat  v^\epsilon  \varphi \,dx'dz+ \int_{ \omega\times Z}\hat v^{\epsilon}_3\partial_{z_3} \hat v^\epsilon\,  \varphi\,dx'dz\right)
\medskip
\\
\displaystyle
-\int_{\omega\times Z}\hat P^\epsilon\, {\rm div}_{x'}\varphi'\,dx'dz-\epsilon^{-\delta}\int_{\omega\times Z}\hat P^\epsilon\,{\rm div}_{z'}\varphi'\,dx'dz\medskip\\
\displaystyle-\epsilon^{-1}\int_{\omega\times Z}\hat P^\epsilon\,\partial_{z_3}\varphi_3\,dx'dz
 =\int_{\omega\times Z} f'\cdot \varphi'\,dx'dz+O_\epsilon.
\end{array}
\end{equation}
 Now, we are ready to pass to the limit in the variational formulation depending on the {\bf TPM}. This will be given in the next subsections.
\subsection{Darcy' law for the {\bf HTPM}}\label{Subsec:HTPM}
Next, we give some compactness results about the behavior of the extended sequences $(\tilde v^\epsilon, \tilde P^\epsilon)$ and the unfolded functions $(\hat v^\epsilon, \hat P^\epsilon)$ by using the a priori estimates given in Lemma \ref{Vel_Estimates_lemma}-{\it (i)}, Theorem \ref{mainthmEstimates}-{\it (i)} and Lemma \ref{estimates_hat}-{\it (i)}, respectively.\\

\begin{lemma} \label{lemma_compactnessH}Consider $\delta> 1$, $\gamma\leq \delta$ and $C(\delta)$ defined by (\ref{C_case_H}).  Then, there exists:
 \begin{itemize}
 \item  a subsequence, still denoted by $(\tilde v^\epsilon, \tilde P^\epsilon)$, chosen from a sequence of $(\tilde v^\epsilon, \tilde P^\epsilon)$ of solutions of (\ref{N-S-d}) (we also denote by $(\hat v^\epsilon, \hat P^\epsilon)$ the subsequence of the corresponding unfolded functions),
 \item $(\tilde v, \tilde P)\in   L^2(\Omega)^3)\times L^{C(\delta)}_0(\omega)$, where  $\tilde v=0$ on $z_3=\{0,1\}$, $\tilde v_3=0$ and $\tilde P$ independent of $z_3$, 
 \item $\hat v\in L^2(\Omega; H^1_{\#}(Z')^3)$, with $\hat v=0$ on $\omega\times T$, $\hat v_3$  independent of $z_3$,  and  satisfying the relation 
 \begin{equation}\label{relation_mean2H}\int_{Z}\hat v(x',z)\,dz=\int_0^1\tilde v(x',z_3)\,dz_3\quad \hbox{with }\int_Z\hat v_3(x',z)\,dz=0,
 \end{equation}
 such that 
\begin{equation}\label{conv_vel_tilde2H}
\epsilon^{\gamma-2\delta}\tilde v^\epsilon \rightharpoonup (\tilde v',0)\hbox{ weakly in }  L^2(\Omega)^3,
\end{equation}
\begin{equation}\label{conv_vel_gorro2H}
\epsilon^{\gamma-2\delta}\hat v^\epsilon\rightharpoonup \hat v\hbox{ weakly in }L^2(\Omega; H^1(Z')^3),\end{equation}
\begin{equation}\label{conv_pres_tilde2H}
\tilde P^\epsilon\to \tilde P \hbox{ strongly in }L^{C(\delta)}(\Omega),
\end{equation}
\begin{equation}\label{conv_pres_hat2H}
\hat P^\epsilon\to \tilde P \hbox{ strongly in }L^{C(\delta)}(\omega\times Z).
\end{equation}
 \end{itemize}
Moreover, $\tilde v$ and $\hat v$ satisfy the following divergence conditions
\begin{equation}\label{divxproperty2H}
{\rm div}_{x'}\left(\int_0^1\tilde v'(x',z_3)\,dz_3\right)=0\  \hbox{ in }\omega,\quad \left(\int_0^1\tilde v'(x',z_3)\,dz_3\right)\cdot n=0\  \hbox{ on }\partial \omega,
\end{equation}
\begin{equation}\label{divyproperty2H}
\begin{array}{c}\displaystyle
 {\rm div}_{z'}\,\hat v'(x',z)=0\  \hbox{ in }\omega\times Z,\medskip\\
\displaystyle  {\rm div}_{x'}\left(\int_{Z}\hat v'(x',z)\,dz\right)=0\  \hbox{ in }\omega, \quad\left(\int_{Z}\hat v'(x',z)\,dz\right)\cdot n=0\  \hbox{ on }\partial\omega.
\end{array}
\end{equation}
\end{lemma}
\begin{proof} The proof of this result is obtained by arguing similarly to  \cite[Section 5]{Anguiano_SuarezGrau} for the case $a_\epsilon\ll \epsilon$ (see also \cite{Anguiano_SuarezGrau2} for more details), by taking into account that in this case $a_\epsilon=\epsilon^\delta$ with $\delta>1$ and that $\gamma$ is included in the estimates of velocity (\ref{estimates_v_H}), so we omit it.

\end{proof}

Following the proof of Theorem \ref{mainthmEstimates}-$(ii)$ and taking into account the estimates for $\hat v^\epsilon$ given in Lemma \ref{estimates_hat}-$(i)$, we prove that the inertial terms given in the variational formulation (\ref{form_var_hat}) vanish when $\epsilon$ tends to zero.\\

\begin{proposition}\label{InertialhatH}
Consider $\delta>1$ and $\gamma\leq \delta$. Then, the inertial term satisfies
\begin{equation}\label{inertialhatHestim}
\begin{array}{l}\displaystyle
\left|-\epsilon^{-\delta}\int_{ \omega\times Z} \hat v^\epsilon\tilde \otimes \hat v^\epsilon:D_{z'}  \varphi\,dx'dz\right.\medskip\\
\displaystyle
\displaystyle \quad \left.+\epsilon^{-1}\left(\int_{ \omega\times Z}\partial_{z_3} \hat v^{\epsilon}_3\hat  v^\epsilon  \varphi \,dx'dz+ \int_{ \omega\times Z}\hat v^{\epsilon}_3\partial_{z_3} \hat v^\epsilon\,  \varphi\,dx'dz\right)\right|\leq C\epsilon^{3\delta-2\gamma}\to 0,
\end{array}
\end{equation}
when $\epsilon$ tends to zero, where $\varphi\in \mathcal{D}(\omega;C^\infty_{\#}(Z)^3)$ and  $(\hat v^\epsilon \tilde \otimes \hat v^\epsilon)_{ij}=\hat v^\epsilon_i \hat v^\epsilon_j$, $i=1,2$, $j=1,2,3$.
\end{proposition}
\begin{proof}
We begin estimating the term in the left-hand side of (\ref{inertialhatHestim}). By using H${\rm \ddot{o}}$lder's inequality and estimates of $\hat v^\epsilon$ given in Lemma \ref{estimates_hat}-$(i)$, we get
\begin{equation}\label{Hestim1}
\begin{array}{rl}\displaystyle
\left|\epsilon^{-\delta}\int_{ \omega\times Z} \hat v^\epsilon\tilde \otimes \hat v^\epsilon:D_{z'}  \varphi\,dx'dz\right|&\displaystyle 
\leq \epsilon^{-\delta}\|\hat v^\epsilon\|^2_{L^2(\omega\times Z)^3}\|D_{z'}\varphi\|_{L^\infty(\omega\times Z)^{3\times 2}}\medskip\\
&\displaystyle\leq C\epsilon^{3\delta-2\gamma}.
\end{array}
\end{equation}
Next, we estimate the rest of the terms of the left-hand side of (\ref{inertialhatHestim}). By using H${\rm \ddot{o}}$lder's inequality and estimates of $\hat v^\epsilon$ given in Lemma \ref{estimates_hat}-$(i)$, we obtain
\begin{equation}\label{Hestim3}
\begin{array}{l}\displaystyle
\left|\epsilon^{-1}\left(\int_{ \omega\times Z}\partial_{z_3} \hat v^{\epsilon}_3\hat  v^\epsilon  \varphi \,dx'dz+ \int_{ \omega\times Z}\hat v^{\epsilon}_3\partial_{z_3} \hat v^\epsilon\,  \varphi\,dx'dz\right)\right| \medskip\\
\displaystyle \leq \displaystyle \epsilon^{-1}\|\hat v^\epsilon\|_{L^2(\omega\times Z)^3}\|\partial_{z_3}\hat v^\epsilon\|_{L^2(\omega\times Z)^3}\|\varphi\|_{L^\infty(\omega\times Z)^3}
\medskip\\
\displaystyle\leq \displaystyle C\epsilon^{3\delta-2\gamma}.
\end{array}
\end{equation}
Since we have assumed $\gamma\leq \delta$, then $3\delta-2\gamma\geq \delta>1$, and so, $\epsilon^{3\delta-2\gamma}\to 0$ when $\epsilon$ tends to zero.

\end{proof}
\noindent Now, by using Lemma \ref{lemma_compactnessH} and Proposition \ref{InertialhatH}, we will give some details of the proof of Theorem \ref{mainthmHTPM}. \\

\noindent {\bf Proof of Theorem \ref{mainthmHTPM}.} Consider $\delta>1$ and $\gamma\leq \delta$. First, we pass to the limit in (\ref{form_var_hat}) to obtain a two pressure homogenized problem. To do this, we choose $\varphi\in \mathcal{D}(\omega; C^\infty_{\#}(Z)^3)$ with ${\rm div}_{z'}\varphi'=0$ in $\omega\times Z$, ${\rm div}_{x'}(\int_{Z}\varphi'\,dz)=0$ in $\omega$ and $\varphi_3$ independent of $y_3$ in (\ref{form_var_hat}). Taking into account that thanks to ${\rm div}_{z'}\varphi'=0$ in $\omega\times Z_f$ and $\varphi_3$ independent of $y_3$, we have 
$$\epsilon^{-\delta}\int_{\omega\times Z}\hat P^\epsilon{\rm div}_{z'}\varphi'\,dx'dz=0\quad\hbox{and}\quad \epsilon^{-1}\int_{\omega\times Z}\partial_{z_3}\hat P^\epsilon\partial_{z_3}\varphi_3\,dx'dz=0.$$
Thus, passing to the limit using convergences given in Lemma \ref{lemma_compactnessH}, taking into account that 
$$\left|\epsilon^{\gamma-2}\eta\int_{ \omega\times Z}\partial_{z_3}\hat v^\epsilon:\partial_{z_3}\varphi\,dx'dz\right|\leq C\epsilon^{\delta-1}\to 0,$$
that the inertial terms tend to zero thanks to (\ref{inertialhatHestim}), that 
$$\int_{\omega\times Z}\hat P^\epsilon\, {\rm div}_{x'}\varphi'\,dx'dz\to \int_{\omega\times Z}\tilde P\, {\rm div}_{x'}\varphi'\,dx'dz,$$
and that $\tilde P$ does not depend on $z$ and ${\rm div}_{x'}(\int_{Z}\varphi'dz)=0$, we obtain
\begin{equation}\label{form_var_hatH1}\begin{array}{l}
\displaystyle
 \eta\int_{ \omega\times Z}D_{z'}\hat v:D_{z'}\varphi\,dx'dz =\int_{\omega\times Z} f'\cdot \varphi'\,dx'dz.
\end{array}
\end{equation}
If we focus in  the third component, taking $\varphi\cdot e_3$ in the previous variational formulation, since the boundary condition $\hat v=0$ on $S$, we deduce that $\hat v_3=0$. Now, we take into account that there is no $z_3$-dependence in the obtained variational formulation. For that, we can consider $\varphi$ independent of $z_3$, which implies that $\hat V(x', z')=\int_0^1\hat v(x',z)\,dz_3$ a.e. $(x',z')\in \omega\times Z'$, with $\hat V_3=0$,  satisfies the same variational formulation with integral in $\omega\times Z'$. By density arguments, see \cite{Anguiano_SuarezGrau},  we can deduce that there exists $\tilde q\in L^2(\omega)/\mathbb{R}$ (which coincides with $\tilde P$) and $\hat q^1(x', z') \in L^2(\omega;L^2_\#( Z'_f))$ such that the variational formulation for $\hat V$ is equivalent to problem 
\begin{equation}\label{twopressureH1}
\left\{\begin{array}{rl}
-\Delta_{z'}\hat V'+\nabla_{z'}\hat q=f'(x')-\nabla_{x'}\tilde P(x') &\hbox{in}\ \omega\times Z'_f,
\medskip
\\
{\rm div}_{z'}\hat V'=0&\hbox{in}\ \omega\times Z'_f,
\medskip
\\
\hat V'=0&\hbox{in}\ \omega\times S',
\medskip
\\
{\rm div}_{x'}\left(\int_{Z'_f}\hat V'(x',z')\,dz'\right)=0&\hbox{in}\ \omega,
\medskip
\\
\left(\int_{Z'_f}\hat V'(x',z')\,dz'\right)\cdot n=0&\hbox{on}\ \partial\omega,
\medskip
\\
\hat V'(x',z'),\ \hat q(x',z')\ Z'-\hbox{periodic}.
\end{array}\right.
\end{equation}
Next, taking into account the local problems (\ref{cell_problem_HTPM}), we eliminate the microscopic variable $z'$ in the effective problem (\ref{twopressureH1}). To do this, by using standard arguments of linearity, we consider the following identification
$$\begin{array}{l}\displaystyle 
\hat V(x',z')={1\over \eta}\sum_{i=1}^2\left(f_i(x')-\partial_{x_i}\tilde P(x')\right)w^i(z'),\medskip\\
\displaystyle \hat q(x',z')=\sum_{i=1}^2\left(f_i(x')-\partial_{x_i}\tilde P(x')\right)\pi^i(z'),
\end{array}$$
with $(w^i, \pi^i)$, $i=1,2$, the unique solution of cell problems (\ref{cell_problem_HTPM}). Thanks to the identities (\ref{relation_mean2H}), we deduce that $\tilde V$ is given by the first line of (\ref{DarcyHTPM}).\\

Finally, the divergence condition with respect to the variable $x'$ together with the expression of $\tilde V$
 gives the second line of (\ref{DarcyHTPM}), which has a unique solution $\tilde P\in H^1(\Omega)\cap L^2_0(\omega)$ since $\mathcal{K}$ is positive definite (see \cite[Theorem 2.1]{Anguiano_SuarezGrau2}). 

\par\hfill$\square$

\subsection{Darcy' law for the {\bf PTPM}}\label{Subsec:PTPM}
Next, we give some compactness results about the behavior of the extended sequences $(\tilde v^\epsilon, \tilde P^\epsilon)$ and the unfolded functions $(\hat v^\epsilon, \hat P^\epsilon)$ by using the a priori estimates given in Lemma \ref{Vel_Estimates_lemma}-{\it (ii)}, Theorem \ref{mainthmEstimates}-{\it (ii)} and Lemma \ref{estimates_hat}-{\it (ii)}, respectively.\\

\begin{lemma} \label{lemma_compactnessP}Consider $\delta= 1$, $\gamma\leq 1$ and $C(\delta)$ defined by (\ref{C_case_H}).  Then, there exists:
 \begin{itemize}
 \item  a subsequence, still denoted by $(\tilde v^\epsilon, \tilde P^\epsilon)$, chosen from a sequence of $(\tilde v^\epsilon, \tilde P^\epsilon)$ of solutions of (\ref{N-S-d}) (we also denote by $(\hat v^\epsilon, \hat P^\epsilon)$ the subsequence of the corresponding unfolded functions),
 \item $(\tilde v, \tilde P)\in   H^1(0,1;L^2(\omega)^3)\times L^{C(1)}_0(\omega)$, where  $\tilde v=0$ on $z_3=\{0,1\}$, $\tilde v_3=0$ and $\tilde P$ independent of $z_3$, 
 \item $\hat v\in L^2(\omega; H^1_{\#}(Z)^3)$, with $\hat v=0$ on $\omega\times S$  and on $z_3=\{0,1\}$,  and  satisfying the relation 
 \begin{equation}\label{relation_mean2P}\int_{Z}\hat v(x',z)\,dz=\int_0^1\tilde v(x',z_3)\,dz_3\quad \hbox{with }\int_Z\hat v_3(x',z)\,dz=0,
 \end{equation}
 such that 
\begin{equation}\label{conv_vel_tilde2P}
\epsilon^{\gamma-2}\tilde v^\epsilon \rightharpoonup (\tilde v',0)\hbox{ weakly in } H^1(0,1; L^2(\omega)^3),
\end{equation}
\begin{equation}\label{conv_vel_gorro2P}
\epsilon^{\gamma-2}\hat v^\epsilon\rightharpoonup \hat v\hbox{ weakly in }L^2(\omega; H^1(Z)^3),\end{equation}
\begin{equation}\label{conv_pres_tilde2P}
\tilde P^\epsilon\to \tilde P \hbox{ strongly in }L^{C(1)}(\Omega),
\end{equation}
\begin{equation}\label{conv_pres_hat2P}
\hat P^\epsilon\to \tilde P \hbox{ strongly in }L^{C(1)}(\omega\times Z).
\end{equation}
 \end{itemize}
Moreover, $\tilde v$ and $\hat v$ satisfy the following divergence conditions
\begin{equation}\label{divxproperty2P}
{\rm div}_{x'}\left(\int_0^1\tilde v'(x',z_3)\,dz_3\right)=0\  \hbox{ in }\omega,\quad \left(\int_0^1\tilde v'(x',z_3)\,dz_3\right)\cdot n=0\  \hbox{ on }\partial \omega,
\end{equation}
\begin{equation}\label{divyproperty2P}
\begin{array}{c}
 {\rm div}_{z}\,\hat v(x',z)=0\  \hbox{ in }\omega\times Z,\medskip\\
\displaystyle  {\rm div}_{x'}\left(\int_{Z}\hat v'(x',z)\,dz\right)=0\  \hbox{ in }\omega, \quad\left(\int_{Z}\hat v'(x',z)\,dz\right)\cdot n=0\  \hbox{ on }\partial\omega.
 \end{array}
\end{equation}
\end{lemma}
\begin{proof} The proof of this result is obtained by arguing similarly to  \cite[Section 5]{Anguiano_SuarezGrau} for the case $a_\epsilon\approx \epsilon$ (see also \cite{Anguiano_SuarezGrau2} for more details), by taking into account that in this case $a_\epsilon=\epsilon$ and $\lambda=1$, and that $\gamma$ is included in the estimates of velocity (\ref{estimates_v_P}), so we omit it.

\end{proof}

\noindent Following the proof of Theorem \ref{mainthmEstimates}-$(ii)$ and taking into account the estimates for $\hat v^\epsilon$ given in Lemma \ref{estimates_hat}-$(ii)$, we have the following result.\\

\begin{proposition}\label{InertialhatP}
Consider $\delta=1$ and $\gamma\leq 1$. Then, the inertial term satisfies
\begin{equation}\label{inertialhatPestim}
\begin{array}{l}\displaystyle
\!\!\!\!\!\left|-\epsilon^{-1}\left(\int_{ \omega\times Z} \!\!\hat v^\epsilon\tilde \otimes \hat v^\epsilon:D_{z'}  \varphi\,dx'dz
\displaystyle +\int_{ \omega\times Z}\!\!\partial_{z_3} \hat v^{\epsilon}_3\hat  v^\epsilon  \varphi \,dx'dz+ \int_{ \omega\times Z}\!\!\hat v^{\epsilon}_3\partial_{z_3} \hat v^\epsilon\,  \varphi\,dx'dz\right)\right|\medskip\\
\displaystyle
\!\!\!\!\!\leq C\epsilon^{3-2\gamma}\to 0,
\end{array}
\end{equation}
when $\epsilon$ tends to zero, where $\varphi\in \mathcal{D}(\omega;C^\infty_{\#}(Z)^3)$ and  $(\hat v^\epsilon \tilde \otimes \hat v^\epsilon)_{ij}=\hat v^\epsilon_i \hat v^\epsilon_j$, $i=1,2$, $j=1,2,3$.
\end{proposition}
\begin{proof} The proof follows the lines of Proposition \ref{InertialhatH} by taking into account that the estimates of $\hat v^\epsilon$ and $D_z \hat v^\epsilon$ given in Lemma \ref{estimates_hat}-$(ii)$ and that $\delta=1$, so we omit it. 

\end{proof}
\noindent {\bf Proof of Theorem \ref{mainthmPTPM}.} Consider $\delta=1$ and $\gamma\leq 1$. 
First, we pass to the limit in (\ref{form_var_hat}) to obtain a two pressure homogenized problem. To do this, we choose $\varphi\in \mathcal{D}(\omega; C^\infty_{\#}(Z)^3)$ with ${\rm div}_{z}\varphi=0$ in $\omega\times Z$, ${\rm div}_{x'}(\int_{Z}\varphi'\,dz)=0$ in $\omega$  in (\ref{form_var_hat}). Taking into account that thanks to ${\rm div}_{z}\varphi=0$ in $\omega\times Z$, we have 
$$\epsilon^{-1}\int_{\omega\times Z}\hat P^\epsilon{\rm div}_{z}\varphi\,dx'dz=0.$$
Thus, passing to the limit using convergences given in Lemma \ref{lemma_compactnessP}, taking into account
that the inertial terms tend to zero thanks to (\ref{inertialhatPestim}), that 
$$\int_{\omega\times Z}\hat P^\epsilon\, {\rm div}_{x'}\varphi'\,dx'dz\to \int_{\omega\times Z}\tilde P\, {\rm div}_{x'}\varphi'\,dx'dz,$$
and that $\tilde P$ does not depend on $z$ and ${\rm div}_{x'}(\int_{Z}\varphi'dz)=0$, we obtain
\begin{equation}\label{form_var_hatP1}\begin{array}{l}
\displaystyle
 \eta\int_{ \omega\times Z}D_{z}\hat v:D_{z}\varphi\,dx'dz =\int_{\omega\times Z} f'\cdot \varphi'\,dx'dz.
\end{array}
\end{equation}
By density arguments, see \cite{Anguiano_SuarezGrau}, we can deduce that there exist $q(x')\in L^2(\omega)/\mathbb{R}$ (which coincides with $\tilde P$) and $\hat q(x',z)\in L^2_\#(\omega\times Z_f)$ such that the variational formulation (\ref{form_var_hatP1}) for $\hat v$ is equivalent to problem 
\begin{equation}\label{twopressureP1}
\left\{\begin{array}{rl}
-\Delta_{z}\hat v+\nabla_{z}\hat q=f'(x')-\nabla_{x'}\tilde P(x') &\hbox{in}\ \omega\times Z_f,
\medskip
\\
{\rm div}_{z}\hat v'=0&\hbox{in}\ \omega\times Z_f,
\medskip
\\
\hat v'=0&\hbox{in}\ \omega\times S,
\medskip
\\
{\rm div}_{x'}\left(\int_{Z_f}\hat v'(x',z)\,dz\right)=0&\hbox{in}\ \omega,
\medskip
\\
\left(\int_{Z_f}\hat v'(x',z)\,dz\right)\cdot n=0&\hbox{on}\ \partial\omega,
\medskip
\\
\hat v(x',z),\ \hat q(x',z)\ Z'-\hbox{periodic}.
\end{array}\right.
\end{equation}
Next, we take into account the local problems (\ref{cell_problem_PTPM}) for $i=1,2,3$ with $e_i=(\delta_{1i}, \delta_{2i},0)^t$ with $\delta_{ji}$ is the Kronecker delta. Observe that $w^3=0$ and $\pi^3$ is constant since for $i=3$ the force term is absent.  Next, we eliminate the microscopic variable  in the effective problem (\ref{twopressureP1}). To do this, we consider the following identification
$$\begin{array}{l}\displaystyle 
\hat v(x',z)={1\over \eta}\sum_{i=1}^2\left(f_i(x')-\partial_{x_i}\tilde P(x')\right)w^i(z'),\medskip\\
\displaystyle \hat q(x',z)=\sum_{i=1}^2\left(f_i(x')-\partial_{x_i}\tilde P(x')\right)\pi^i(z'),
\end{array}$$
and thanks to the identities (\ref{relation_mean2P}), we deduce that $\tilde V$ is given by 
$$\tilde V(x')={1\over \eta}\mathcal{K}_P\left(f'(x')-\nabla_{x'}\tilde P(x')\right), \hbox{ with}\quad \mathcal{K}_P=\int_{Z_f}w^i_j(z)\,dz\in \mathbb{R}^{3\times 3},\quad  i,j=1,2,3.$$
Observe that $(\mathcal{K}_P)_{i3}=(\mathcal{K}_P)_{3,j}=0$, $i,j=1,2,3$, so that $\tilde V$  can be rewritten to have (\ref{DarcyPTPM}) with $\mathcal{K}_P$ given by (\ref{JPTPM}).\\

Finally, the divergence condition with respect to the variable $x'$ together with the expression of $\tilde V$
 gives the second line of (\ref{DarcyPTPM}), which has a unique solution $\tilde P\in H^1(\Omega)\cap L^2_0(\omega)$ since $\mathcal{K}_{P}$ is positive definite (see \cite[Theorem 2.1]{Anguiano_SuarezGrau2}).
\par\hfill$\square$

\subsection{Darcy' law for the {\bf VTPM}}\label{Subsec:VTPM}
We give some compactness results about the behavior of the extended sequences $(\tilde v^\epsilon, \tilde P^\epsilon)$ and the unfolded functions $(\hat v^\epsilon, \hat P^\epsilon)$ by using the a priori estimates given in Lemma \ref{Vel_Estimates_lemma}-{\it (iii)}, Theorem \ref{mainthmEstimates}-{\it (iii)} and Lemma \ref{estimates_hat}-{\it (iii)}, respectively.\\

\begin{lemma} \label{lemma_compactnessV}Consider $\delta\in (0, 1)$, $\gamma\leq 1$ and $C(\delta)$ defined by (\ref{C_case_H}).  Then, there exists:
 \begin{itemize}
 \item  a subsequence, still denoted by $(\tilde v^\epsilon, \tilde P^\epsilon)$, chosen from a sequence of $(\tilde v^\epsilon, \tilde P^\epsilon)$ of solutions of (\ref{N-S-d}) (we also denote by $(\hat v^\epsilon, \hat P^\epsilon)$ the subsequence of the corresponding unfolded functions),
 \item $(\tilde v, \tilde P)\in   H^1(0,1;L^2(\omega)^3)\times L^{C(1)}_0(\omega)$, where  $\tilde v=0$ on $z_3=\{0,1\}$, $\tilde v_3=0$ and $\tilde P$ independent of $z_3$, 
 \item $\hat v\in H^1(0,1;L^2_{\#}(\omega\times Z')^3)$, with $\hat v=0$ on $\omega\times S$  and on $z_3=\{0,1\}$,  $\hat v_3$ independent of $z_3$, and  satisfying the relation 
 \begin{equation}\label{relation_mean2V}\int_{Z}\hat v(x',z)\,dz=\int_0^1\tilde v(x',z_3)\,dz_3\quad \hbox{with }\int_Z\hat v_3(x',z)\,dz=0,
 \end{equation}
 such that 
\begin{equation}\label{conv_vel_tilde2V}
\epsilon^{\gamma-2}\tilde v^\epsilon \rightharpoonup (\tilde v',0)\hbox{   in } H^1(0,1; L^2(\omega)^3),
\end{equation}
\begin{equation}\label{conv_vel_gorro2V}
\epsilon^{\gamma-2}\hat v^\epsilon\rightharpoonup \hat v\hbox{   in }H^1(0,1;L^2(\omega\times Z')^3),\end{equation}
\begin{equation}\label{conv_vel_gorro2V2}
\epsilon^{\gamma-2\delta}D_{z'}\hat v\rightharpoonup 0\hbox{   in }L^2(\omega\times Z')^{3\times 2},
\end{equation}
\begin{equation}\label{conv_pres_tilde2V}
\tilde P^\epsilon\to \tilde P \hbox{   in }L^{C(1)}(\Omega),
\end{equation}
\begin{equation}\label{conv_pres_hat2V}
\hat P^\epsilon\to \tilde P \hbox{   in }L^{C(1)}(\omega\times Z).
\end{equation}
 \end{itemize}
Moreover, $\tilde v$ and $\hat v$ satisfy the following divergence conditions
\begin{equation}\label{divxproperty2V}
{\rm div}_{x'}\left(\int_0^1\tilde v'(x',z_3)\,dz_3\right)=0\  \hbox{ in }\omega,\quad \left(\int_0^1\tilde v'(x',z_3)\,dz_3\right)\cdot n=0\  \hbox{ on }\partial \omega,
\end{equation}
\begin{equation}\label{divyproperty2V}
\begin{array}{c}\displaystyle
 {\rm div}_{z'}\,\hat v'(x',z)=0\  \hbox{ in }\omega\times Z,\medskip\\
\displaystyle   {\rm div}_{x'}\left(\int_{Z}\hat v'(x',z)\,dz\right)=0\  \hbox{ in }\omega, \quad\left(\int_{Z}\hat v'(x',z)\,dz\right)\cdot n=0\  \hbox{ on }\partial\omega.
\end{array}
\end{equation}
\end{lemma}
\begin{proof} The proof of this result is obtained by arguing similarly to  \cite[Section 5]{Anguiano_SuarezGrau} for the case $a_\epsilon\gg \epsilon$ (see also \cite{Anguiano_SuarezGrau2} for more details), by taking into account that in this case $a_\epsilon=\epsilon^\delta$ with $\delta\in(0,1)$ and that and that $\gamma$ is included in the estimates of velocity (\ref{estimates_v_V}), so we omit it.

\end{proof}

\noindent Following the proof of Theorem \ref{mainthmEstimates}-$(iii)$ and taking into account the estimates for $\hat v^\epsilon$ given in Lemma \ref{estimates_hat}-$(iii)$, we have the following result.\\

\begin{proposition}\label{InertialhatP}
Consider $\delta\in (0,1)$ and $\gamma\leq 1$. Then, the inertial term satisfies
\begin{equation}\label{inertialhatVestim}
\begin{array}{l}\displaystyle 
\left|-\epsilon^{-\delta}\int_{ \omega\times Z} \hat v^\epsilon\tilde \otimes \hat v^\epsilon:D_{z'}  \varphi\,dx'dz
\displaystyle \right. \medskip\\
\displaystyle\quad \left.+\epsilon^{-1}\left(\int_{ \omega\times Z}\partial_{z_3} \hat v^{\epsilon}_3\hat  v^\epsilon  \varphi \,dx'dz+ \int_{ \omega\times Z}\hat v^{\epsilon}_3\partial_{z_3} \hat v^\epsilon\,  \varphi\,dx'dz\right)\right|\leq  C\epsilon^{3-2\gamma}\to 0,
\end{array}
\end{equation}
when $\epsilon$ tends to zero, where $\varphi\in \mathcal{D}(\omega;C^\infty_{\#}(Z)^3)$ and  $(\hat v^\epsilon \tilde \otimes \hat v^\epsilon)_{ij}=\hat v^\epsilon_i \hat v^\epsilon_j$, $i=1,2$, $j=1,2,3$.
\end{proposition}
\begin{proof}
The proof follows the lines of Proposition \ref{InertialhatH} taking into account that $\delta \in (0,1)$ and estimates of $\hat v^\epsilon$ and $\partial_{z_3}\hat v^\epsilon$ given in Lemma \ref{estimates_hat}-$(iii)$, so we omit it. 

\end{proof}

\noindent {\bf Proof of Theorem \ref{mainthmVTPM}.} Consider $\delta\in (0,1)$ and $\gamma\leq 1$. 
First, we pass to the limit in (\ref{form_var_hat}) to obtain a two pressure homogenized problem. To do this, we choose $\varphi\in \mathcal{D}(\omega; C^\infty_{\#}(Z)^3)$ with ${\rm div}_{z'}\varphi'=0$ in $\omega\times Z$, ${\rm div}_{x'}(\int_{Z}\varphi'\,dz)=0$ in $\omega$ and $\varphi_3$ independent of $z_3$  in (\ref{form_var_hat}). Taking into account that thanks to ${\rm div}_{z'}\varphi'=0$ in $\omega\times Z$ and $\varphi_3$ independent of $z_3$, we have  that
$$\epsilon^{-\delta}\int_{\omega\times Z}\hat P^\epsilon{\rm div}_{z'}\varphi'\,dx'dz=0\quad \hbox{and}\quad\epsilon^{-1}\int_{\omega\times Z}\hat P^\epsilon\partial_{z_3}\varphi_3\,dx'dz=0.$$
Thus, passing to the limit using convergences given in Lemma \ref{lemma_compactnessV}, taking into account
that 
$$\epsilon^{\gamma-2\delta}\eta\int_{ \omega\times Z}D_{z'}\hat v^\epsilon:D_{z'}\varphi\,dx'dz\to 0,$$
that the inertial terms tend to zero thanks to (\ref{inertialhatVestim}), that 
$$\int_{\omega\times Z}\hat P^\epsilon\, {\rm div}_{x'}\varphi'\,dx'dz\to \int_{\omega\times Z}\tilde P\, {\rm div}_{x'}\varphi'\,dx'dz,$$
and that $\tilde P$ does not depend on $z$ and ${\rm div}_{x'}(\int_{Z}\varphi'dz)=0$, we obtain
\begin{equation}\label{form_var_hatV1}\begin{array}{l}
\displaystyle
 \eta\int_{ \omega\times Z}\partial_{z_3}\hat v':\partial_{z_3}\varphi\,dx'dz =\int_{\omega\times Z} f'\cdot \varphi'\,dx'dz.
\end{array}
\end{equation}
 By density arguments, see \cite{Anguiano_SuarezGrau}, we can deduce that there exist $q(x')\in L^2_0(\omega)$ (which coincides with $\tilde P$) and $\hat q(x',z')\in L^2_\#(\omega\times Z')$, such that the variational formulation (\ref{form_var_hatV1}) for $\hat v$ is equivalent to problem 
\begin{equation}\label{twopressureV1}
\left\{\begin{array}{rl}
-\eta\partial^2_{z_3}\hat v'+\nabla_{z'}\hat q=f'(x')-\nabla_{x'}\tilde P(x') &\hbox{in}\ \omega\times Z_f,
\medskip
\\
{\rm div}_{z'}\hat v'=0&\hbox{in}\ \omega\times Z_f,
\medskip
\\
\hat v'=0&\hbox{in}\ \omega\times S,
\medskip
\\
\hat v'=0&\hbox{in}\ z_3=\{0,1\},
\medskip
\\
{\rm div}_{x'}\left(\int_{Z_f}\hat v'(x',z)\,dz\right)=0&\hbox{in}\ \omega,
\medskip
\\
\left(\int_{Z_f}\hat v'(x',z)\,dz\right)\cdot n=0&\hbox{on}\ \partial\omega,
\medskip
\\
\hat v(x',z),\ \hat q(x',z')\ Z'-\hbox{periodic}.
\end{array}\right.
\end{equation}
We just remark that since $\hat v_3$ is independent of $z_3$ and due to the boundary conditions $\hat v_3=0$ on $z_3=\{0,1\}$, then $\hat v_3=0$.\\

\noindent Next, we consider  $(w^i, \pi^i)\in H^1_{\#}(Z_f)^2\times H^1_{\#}(Z'_f)$, $i=1,2$, which is the unique solution of 
\begin{equation}\label{cell_problem_VTPM}
\left\{\begin{array}{rl}
-\partial^2_{z_3}w^i + \nabla_{z'}\pi^i=-e_i &\hbox{in}\  Z_f,
\medskip\\
{\rm div}_{z'} w^i=0&\hbox{in}\  Z_f,
\medskip\\
w^i=0&\hbox{in}\  S,
\medskip\\
w^i=0&\hbox{on}\  z_3=\{0,1\}.
\end{array}
\right.
\end{equation}
By using now the identities
$$\begin{array}{l}\displaystyle\hat v'(x',z)=-{1\over \eta}\sum_{i=1}^2\left(f_i(x')-\partial_{x_i}\tilde P(x')\right)w^i(z),\medskip\\
\displaystyle\hat q(x',z)=-\sum_{i=1}^2\left(f_i(x')-\partial_{x_i}\tilde P(x')\right)\pi^i(z'),
\end{array}$$
we can obtain that the average velocity $\tilde V$ is given by
$$\tilde V(x')=-{1\over \eta}\mathcal{K}_V(f'(x')-\nabla_{x'}\tilde P(x')),\quad \tilde V_3(x')=0,$$
with $\mathcal{K}$ is given by 
$$(\mathcal{K}_V)_{ij}=\int_{Z_f}w^i_j\,dz,\quad i,j=1,2.$$ Then, by the divergence condition in the variable $x'$, we get the Darcy equation
$${\rm div}_{x'} \tilde V'(x')=0\ \hbox{in}\ \omega, \quad \tilde V'(x')\cdot n=0\ \hbox{on}\ \partial\omega.$$
\noindent However, we can obtain a   simpler formula for $\mathcal{K}_V$. We observe that (\ref{cell_problem_VTPM})  can be viewed as a system of ordinary differential equations with constant coefficient, with respect to the variable $z_3$ and unknowns functions $z_3\mapsto w^i_1, w^i_2$, $i=1,2$, where $z'$ is a parameter, $z'\in Z'$. Thus, we can give an explicit expressions for $w^i$ given in terms of $\pi^i$ as follows
$$w^i(z)={1\over 2}(z_3^2-z_3)\left(\nabla_{z'}\pi^i(z')+e_i\right),\quad i=1,2,
$$
and imposing the divergence condition in (\ref{cell_problem_VTPM}), then $\pi^i$, $i=1,2$, satisfies the Hele-Shaw problem 
$$-\Delta_{z'} \pi^i(z')=0\quad\hbox{in } Z'_f,\quad \left(\nabla_{z'}\pi^i(z')+e_i\right)\cdot n=0\quad\hbox{ on }\partial S'.$$
\noindent Integrating $w^i$ with respect to $z_3$, we get
$$\int_0^1w^i_j(z)\,dz=-{1\over 12}\left(\nabla_{z'}\pi^i(z')+e_i\right)e_j,\quad i,j=1,2,$$
so we get that the averaged velocity $\tilde V$ is given by (\ref{DarcyVTPM}).  Moreover, it holds that there is a unique solution $\tilde P\in H^1(\Omega)\cap L^2_0(\omega)$ of (\ref{DarcyVTPM}) because $\mathcal{K}_V$ is definite positive, see \cite[Theorem 2.1]{Anguiano_SuarezGrau2}, which finishes the proof.

\par\hfill$\square$

\subsection*{Conflict of interest.} The authors confirm that there is no conflict of interest to report.

\end{document}